\documentclass{tac}
\usepackage[all,ps,arc]{xy}  
\usepackage{amsmath,amssymb,latexsym}
\usepackage{pdfsync}
\usepackage{makeidx}
\usepackage[english]{babel}
\usepackage[applemac]{inputenc}
\usepackage{graphicx}
\usepackage[x11names]{xcolor}
\usepackage{hyperref}

\newtheorem{Theorem}{Theorem}[section]
\newtheorem{Lemma}[Theorem]{Lemma}
\newtheorem{Proposition}[Theorem]{Proposition}

\newtheorem{Corollary}[Theorem]{Corollary}

\newtheorem{Text}[Theorem]{}
\newcommand{\bA}{{\mathbb A}}
\newcommand{\bB}{{\mathbb B}}
\newcommand{\bC}{{\mathbb C}}
\newcommand{\bD}{{\mathbb D}}
\newcommand{\bK}{{\mathbb K}}
\newcommand{\bP}{{\mathbb P}}
\newcommand{\bR}{{\mathbb R}}
\newcommand{\bV}{{\mathbb V}}
\newcommand{\bX}{{\mathbb X}}
\newcommand{\cA}{{\mathcal A}}
\newcommand{\cB}{{\mathcal B}}
\newcommand{\cH}{{\mathcal H}}
\newcommand{\cN}{{\mathcal N}}
\newcommand\Arr{\mathbf{Arr}}
\newcommand\Grpd{\mathbf{Grpd}}
\newcommand\Ker{\mathrm{Ker}}
\newcommand\bKer{\mathbb K\mathrm{er}}
\newcommand\K{\mathrm{K}}
\newcommand\id{\mathrm{id}}
\newcommand\Id{\mathrm{Id}}
\newcounter{functor}

\title{On fibrations between internal groupoids}
\author{P.-A. Jacqmin, S. Mantovani, G. Metere, E.M. Vitale}
\amsclass{18A22, 18A30, 18B40, 18G55}
\keywords{internal groupoid, fibration, strong h-pullback, protomodular category}
\address
{Institut de recherche en math\'ematique et physique, Universit\'e catholique de Louvain\\
Chemin du Cyclotron 2, B 1348 Louvain-la-Neuve, Belgique.\\
Dipartimento di matematica, Universit\`a degli studi di Milano\\
Via C. Saldini 50, 20133 Milano, Italia.\\
Dipartimento di matematica e informatica, Universit\`a degli studi di Palermo\\
Via Archirafi 34, 90123 Palermo, Italia.}
\copyrightyear{2017}
\thanks{Financial support from FNRS grant 1.A741.14 is gratefully acknowledged by the first author. The second and the third author acknowledge the financial support of the I.N.D.A.M.\ Gruppo Nazionale per le Strutture Algebriche, Geometriche e le loro Applicazioni.}
\eaddress{pierre-alain.jacqmin@uclouvain.be,\\ sandra.mantovani@unimi.it,\\ giuseppe.metere@unipa.it,\\ enrico.vitale@uclouvain.be}

\begin{document}

\maketitle


\begin{abstract}
In order to deduce the internal version of the Brown exact sequence from the internal version of the Gabriel-Zisman exact 
sequence, we characterize fibrations and $\ast$-fibrations in the 2-category of internal groupoids in terms of the comparison 
functor from certain pullbacks to the corresponding strong homotopy pullbacks. A similar analysis in the category 
of arrows allows us to give a characterization of protomodular categories using strong homotopy kernels.
\end{abstract}

\tableofcontents

\section{Introduction}

In~\cite{GZ}, Gabriel and Zisman constructed a $\pi_0$-$\pi_1$ exact sequence starting from a functor $F \colon \bA \to \bB$
between pointed groupoids. The Gabriel-Zisman sequence involves, as first and fourth point, $\pi_1$ and $\pi_0$ of the strong 
homotopy kernel of $F$. A similar $\pi_0$-$\pi_1$ exact sequence is obtained in~\cite{RB} by Brown, but now the functor 
$F$ is assumed to be a fibration, and the sequence involves, instead of the strong homotopy kernel of $F$, its (strict) kernel.
Both exact sequences have been generalized in~\cite{SnailMMVShort}, replacing pointed functors by functors between groupoids internal to a pointed regular category $\cA$ with reflexive coequalizers. In order to deduce the internal version of the Brown sequence from the internal version of the Gabriel-Zisman sequence, one needs the fact that, if the internal functor $F$ is a fibration, then the comparison functor $J \colon \bKer(F) \to \bK(F)$ from the kernel to the strong homotopy kernel is a weak equivalence (so that the induced arrows $\pi_0(F)$ and $\pi_1(F)$ are isomorphisms).

The first aim of this note is therefore to prove that, if $F$ is a fibration, then $J$ is a weak equivalence. Since the converse 
implication is not true, we work out a more complete analysis of the situation getting the following results. In Section~\ref{section strong pullbacks}, we review
some basic facts on homotopy pullbacks in a 2-category with invertible 2-cells, in order to conclude that for a category $\cA$ with pullbacks, $\Grpd(\cA)$, the 
2-category of internal groupoids, has strong pullbacks and strong homotopy pullbacks. In Section~\ref{section fibrations}, we assume the base category
$\cA$ to be regular and we prove that $F$ is a fibration if and only if the comparison functor from a suitable pullback to the corresponding 
strong homotopy pullback is a weak equivalence. Here, the pullback and the strong homotopy pullback are those of $F$ along the embedding 
of $B_0$, the object of objects of $\bB$, into $\bB$. In Section~\ref{section *-fibrations}, we assume that $\cA$ is regular and pointed, and we characterize 
$\ast$-fibrations, that is, those functors $F$ such that the comparison from the kernel to the strong homotopy kernel is a weak equivalence.
Thanks to the regularity of $\cA$, any fibration is a $\ast$-fibration, so that we get as a corollary the result needed to compare
Gabriel-Zisman and Brown sequences.

The normalized version of Brown and Gabriel-Zisman sequences are the snake and snail sequences, studied in the context of regular
protomodular categories in~\cite{DB, SnailEV, SnailZJEV}. In this context, the condition of being a fibration is replaced by the 
condition that a certain arrow is a regular epimorphism. In order to have an analysis as complete as possible of the situation in the 
normalized context, in Section~\ref{section normalized fibrations} we compare the notions of fibration, $\ast$-fibration and weak equivalence in $\Grpd(\cA)$ 
with suitable notions of fibration, $\ast$-fibration and weak equivalence in $\Arr(\cA)$, the category with null-homotopies 
whose objects are arrows in $\cA$. We show that, if $\cA$ is pointed, regular and protomodular, the normalization process 
$\cN \colon \Grpd(\cA) \to \Arr(\cA)$ preserves and reflects fibrations, $\ast$-fibrations and weak equivalences.
(We recall that in the more specific case with $\cA$ semi-abelian, the normalization process yields an equivalence between $\Grpd(\cA)$ and the category $\mathbf{XMod}(\cA)$ of internal crossed-modules, \cite{Janelidze2003}).

Moreover, in studying the relation between fibrations and $\ast$-fibrations in $\Arr(\cA)$, with $\cA$  pointed and regular, we found the unexpected result that  
the implication ``fibration $\Rightarrow$ $\ast$-fibration''
is in fact equivalent to the condition that $\cA$ is protomodular.

\hfill

Note that in this paper, the composition of two arrows
$$\xymatrix{ \ar[r]^-{f} & \ar[r]^-{g} &}$$
will be denoted by $f \cdot g$.

\section{Strong pullbacks and strong h-pullbacks}\label{section strong pullbacks}

\begin{Text}{\rm 
Let  $\underline\cB$ be a 2-category with invertible 2-cells, and $\cB$ its underlying category.
We adopt the following terminology:\\
1. A 1-cell $F \colon \bA \to \bB$ in $\underline\cB$ is {\em fully faithful} if, for any $\bX$ in $\underline\cB$, the induced functor 
$$- \cdot F \colon \underline\cB(\bX,\bA) \to \underline\cB(\bX,\bB)$$
is fully faithful in the usual sense.\\
2. Consider 1-cells $F \colon \bA \to \bB$ and $G \colon \bC \to \bB$ in $\underline\cB$.
A {\em strong homotopy pullback} (strong h-pullback, for short) of $F$ and $G$ is a diagram of the form
$$\xymatrix{\bP \ar[r]^{G'} \ar[d]_{F'} & \bA \ar[d]^{F} \\
\bC \ar@{}[ru]^{\varphi}|{\Rightarrow} \ar[r]_{G} & \bB}$$
satisfying the following universal property:
\begin{enumerate}
\item[(a)] For any diagram of the form 
$$\xymatrix{\bX \ar[r]^{H} \ar[d]_{K} & \bA \ar[d]^{F} \\
\bC \ar@{}[ru]^{\mu}|{\Rightarrow} \ar[r]_{G} & \bB}$$
there exists a unique 1-cell $T \colon \bX \to \bP$ such that $T \cdot G' = H, T \cdot F' = K$ and $T \cdot \varphi = \mu$.

\item[(b)] Given 1-cells $L, M \colon \bX \rightrightarrows \bP$ and 2-cells $\alpha \colon L \cdot F' \Rightarrow  M \cdot F'$ and
$\beta \colon L \cdot G' \Rightarrow M \cdot G'$, if
$$\xymatrix{L \cdot F' \cdot G \ar@{=>}[r]^-{\alpha \cdot G} \ar@{=>}[d]_{L \cdot \varphi} & M \cdot F' \cdot G \ar@{=>}[d]^{M \cdot \varphi} \\
L \cdot G' \cdot F \ar@{=>}[r]_-{\beta \cdot F} & M \cdot G' \cdot F}$$
commutes, then there exists a unique 2-cell $\mu \colon L \Rightarrow M$ such that $\mu \cdot F' = \alpha$ and $\mu \cdot G' = \beta$.
\end{enumerate}

3. We say that a pullback $\bC \times_{G,F} \bA$ of $F \colon \bA \rightarrow \bB$ and $G \colon \bC \rightarrow \bB$ in the category 
$\cB$ is {\em strong} (in $\underline\cB$) if
$$\xymatrix{\bC \times_{G,F} \bA \ar[r]^-{\widehat{G}} \ar[d]_-{\widehat{F}} & \bA \ar[d]^-{F} \\
\bC \ar@{}[ru]^{\id}|{\Rightarrow} \ar[r]_-{G} & \bB}$$
satisfies condition (b) above.
}\end{Text}

\begin{Text}\label{TextStrongHpb}{\rm
Another way to express the universal property of the strong h-pullback is first to fix an object $\bX$ in $\underline\cB$ and to construct 
the comma-square of groupoids (which is precisely the strong h-pullback in the 2-category of groupoids).
$$\xymatrix{(- \cdot G \downarrow - \cdot F) \ar[rr] \ar[d] & \ar@{}[d]|{\cong} & \underline\cB(\bX,\bA) \ar[d]^{- \cdot F} \\
\underline\cB(\bX,\bC) \ar[rr]_-{- \cdot G} & & \underline\cB(\bX,\bB)}$$
Then the universal property of the strong h-pullback means that, for any $\bX$, the canonical comparison functor
$$\underline\cB(\bX,\bP) \to (- \cdot G \downarrow - \cdot F)$$
is bijective on objects (condition a) and fully faithful (condition b), that is, it is an isomorphism of categories. 
This makes evident that a strong h-pullback is determined by its universal property up to isomorphism. \\
A weaker universal property consists in asking that the canonical comparison functors $\underline\cB(\bX,\bP) \to (- \cdot G \downarrow - \cdot F)$
are equivalences of groupoids. In this way one gets what is sometimes called a {\em bipullback}, which is determined only up to equivalence.\\
Intermediate situations between strong homotopy pullbacks and bipullbacks are considered in the literature. For example, in~\cite{GR} the
comparison functors are required to be bijective on objects but not fully faithful (the name of h-pullback is used in this case), 
and in~\cite{GZ} the comparison functors are required to be surjective on objects and full.
}\end{Text}

\begin{Text}\label{TextVectB}{\rm
Among strong h-pullbacks, the following one plays a special role.
$$\xymatrix{\vec{\bB} \ar[r]^{\gamma} \ar[d]_{\delta} & \bB \ar[d]^{\Id} \\
\bB \ar@{}[ru]^{\beta}|{\Rightarrow} \ar[r]_{\Id} & \bB}$$
Indeed, if for a category $\cA$ we denote by $\Arr(\cA)$ the category having arrows of $\cA$ as objects and 
commutative squares as arrows, then the universal property of $\vec\bB$ gives an isomorphism of categories
$$\underline\cB(\bX,\vec\bB) \to \Arr(\underline\cB(\bX,\bB))$$
so that to give a 1-cell $\bX \to \vec\bB$ is the same as giving a 2-cell 
$\xymatrix{\bX \ar@{}[r]|{\Downarrow} \ar@/^0,7pc/[r] \ar@/_0,7pc/[r] & \bB}$.
We refer the reader to Section~\ref{section normalized fibrations} for a more detailed treatment of the category $\Arr(\cA)$.
}\end{Text}

\begin{Text}\label{TextPBStrongHpb}{\rm 
Pasting together two strong h-pullbacks, in general one does not get a strong h-pullback. The main interest of the notion of strong pullback 
relies on the following fact: given a diagram in $\underline\cB$ of the form
$$\xymatrix{\bP' \ar[r]^-{\widehat H} \ar[d]_{\widehat{F'}} & \bP \ar[rr]^{G'} \ar[d]^{F'} & & \bA \ar[d]^{F} \\
\bD \ar[r]_-{H} \ar@{}[ru]^{\id}|{\Rightarrow} & \bC \ar@{}[rru]^{\varphi}|{\Rightarrow} \ar[rr]_{G} & & \bB}$$
if  the right-hand part is a strong h-pullback, then the total diagram is a strong h-pullback if and only if the left-hand part is a strong pullback.\\
This fact has an interesting consequence on the existence of strong h-pullbacks: assume that $\underline\cB$ has strong pullbacks
and that the strong h-pullback 
$$\xymatrix{\vec{\bB} \ar[r]^{\gamma} \ar[d]_{\delta} & \bB \ar[d]^{\Id} \\
\bB \ar@{}[ru]^{\beta}|{\Rightarrow} \ar[r]_{\Id} & \bB}$$
exists in $\underline\cB$ for any object $\bB$. Then, for any pair of 1-cells $F \colon \bA \to \bB, G \colon \bC \to \bB$, 
a strong h-pullback of $F$ and $G$ exists and can be obtained by the following limit of solid arrows in $\cB$
$$\xymatrix{& & \bP \ar[lld]_{F'} \ar[d]^{\phi} \ar[rrd]^{G'} \\
\bC \ar[rd]_{G} & & \vec\bB \ar[ld]_{\delta} \ar[rd]^{\gamma} & & \bA \ar[ld]^{F} \\
& \bB\ar@{-->}[dr]_{id} \ar@{}[rr]^{\beta}|{\Rightarrow}& & \bB\ar@{-->}[dl]^{id}\\
&&\bB}$$
together with $\varphi = \phi \cdot \beta \colon F' \cdot G = \phi \cdot \delta \Rightarrow \phi \cdot \gamma = G' \cdot F$.
}\end{Text}

\begin{Text}\label{TextStability}{\rm
Later, we will use the following facts on strong pullbacks and strong h-pullbacks.\\
1. If the pullback in $\cB$ and the strong h-pullback in $\underline\cB$ of $F$ and $G$ exist,
$$\xymatrix{\bC \times_{G,F} \bA \ar[rd]^-{T} \ar@/^1,5pc/[rrd]^-{\widehat{G}} \ar@/_1,5pc/[rdd]_-{\widehat{F}} \\
& \bP \ar[r]^{G'} \ar[d]_{F'} & \bA \ar[d]^{F} \\
& \bC \ar@{}[ru]^{\varphi}|{\Rightarrow} \ar[r]_{G} & \bB}$$
then the pullback $\bC \times_{G,F} \bA$ is strong if and only if the canonical comparison $T$ is fully faithful.\\
2. Consider the pullback in $\cB$ and the strong h-pullback in $\underline\cB$ of $F$ and $G$.
$$\xymatrix{\bC \times_{G,F} \bA \ar[d]_{\widehat F} \ar[r]^-{\widehat G} & \bA \ar[d]^{F} \\
\bC \ar[r]_-{G} & \bB} \; \; \; \; \; \; \; \; \; \; 
\xymatrix{\bP \ar[r]^{G'} \ar[d]_{F'} & \bA \ar[d]^{F} \\
\bC \ar@{}[ru]^{\varphi}|{\Rightarrow} \ar[r]_{G} & \bB}$$
Clearly, if $F$ is fully faithful, then $F'$ is fully faithful. Moreover, if the pullback $\bC \times_{G,F} \bA$ is strong and if 
$F$ is fully faithful, then $\widehat F$ is fully faithful.
}\end{Text}

\begin{Text}\label{TextNotGrpd}{\rm
Now we specialize the previous discussion taking as $\underline\cB$ the 2-category $\Grpd(\cA)$ of groupoids, functors
and natural transformations internal to a category $\cA$ with pullbacks.
The notation for a groupoid $\bB$ in $\cA$ is 
$$\bB=(\xymatrix{B_1 \times_{c,d}B_1 \ar[r]^-{m} & B_1 \ar@<0.8ex>[r]^{d} \ar@<-0.8ex>[r]_{c} & B_0 \ar[l]|{e}} \;,\;
\xymatrix{B_1 \ar[r]^{i} & B_1})$$
where
$$\xymatrix{B_1 \times_{c,d}B_1 \ar[r]^-{\pi_2} \ar[d]_{\pi_1} & B_1 \ar[d]^{d} \\
B_1 \ar[r]_{c} & B_0}$$
is a pullback. The notation for a natural transformation $\alpha \colon F \Rightarrow G \colon \bA \rightrightarrows \bB$ is
$$\xymatrix{A_1 \ar@<0,5ex>[rr]^-{F_1} \ar@<-0,5ex>[rr]_-{G_1} \ar@<-0,5ex>[d]_{d} \ar@<0,5ex>[d]^{c} & & 
B_1 \ar@<-0,5ex>[d]_{d} \ar@<0,5ex>[d]^{c} \\
A_0 \ar@<0,5ex>[rr]^-{F_0} \ar@<-0,5ex>[rr]_-{G_0} \ar[rru]^<<<<<<<<{\alpha} & & B_0}$$
}\end{Text}

\begin{Text}\label{TextStrongHpbGrpd}{\rm
Following~\ref{TextPBStrongHpb}, to prove that $\Grpd(\cA)$ has strong h-pullbacks we need two ingredients. The first one is easy,
the second one is the standard construction of the groupoid of `arrows', and we recall it from~\cite{FractEV} or~\cite{DR}.\\
1. Since pullbacks in $\Grpd(\cA)$ are constructed level-wise, it is straightforward to check that they are strong.\\
2. For every internal groupoid $\bB$, the strong h-pullback
$$\xymatrix{\vec{\bB} \ar[r]^{\gamma} \ar[d]_{\delta} & \bB \ar[d]^{\Id} \\
\bB \ar@{}[ru]^{\beta}|{\Rightarrow} \ar[r]_{\Id} & \bB}$$
exists, and it can be described as follows: 
$$\vec\bB=(\xymatrix{\vec B_1 \times_{\vec c,\vec d}\vec B_1 \ar[r]^-{\vec m} & 
\vec B_1 \ar@<0.9ex>[r]^{\vec d} \ar@<-0.9ex>[r]_{\vec c} & B_1 \ar[l]|{\vec e}} \;,\;
\xymatrix{\vec B_1 \ar[r]^{\vec i} & \vec B_1})$$
where
$$\xymatrix{\vec B_1 \ar[r]^-{m_2} \ar[d]_{m_1} & B_1 \times_{c,d}B_1 \ar[d]^{m} \\
B_1 \times_{c,d}B_1 \ar[r]_-{m} & B_1}$$
is a pullback,  and $\vec d$, $\vec c$ and $\vec e$ are defined by
$$\vec d \colon \xymatrix{\vec B_1 \ar[r]^-{m_1} & B_1 \times_{c,d}B_1 \ar[r]^-{\pi_1} & B_1} \;\;\;\;\;
\vec c \colon \xymatrix{\vec B_1 \ar[r]^-{m_2} & B_1 \times_{c,d}B_1 \ar[r]^-{\pi_2} & B_1}$$
$$\xymatrix{B_1 \ar@/^1pc/[rrd]^{\langle d \cdot e, \id \rangle} \ar[rd]^{\vec e} \ar@/_1pc/[rdd]_{\langle \id, c \cdot e \rangle} \\ 
& \vec B_1 \ar[r]^-{m_2} \ar[d]_{m_1}  & B_1 \times_{c,d} B_1 \ar[d]^{m} \\
& B_1 \times_{c,d} B_1 \ar[r]_-{m} & B_1}$$
The groupoid $\vec\bB$ is equipped with two functors $\delta \colon \vec\bB \to \bB$ and $\gamma \colon \vec\bB \to \bB$ given by
$$\xymatrix{\vec B_1 \ar[rr]^-{\delta_1=m_2 \cdot \pi_1} \ar@<-0.5ex>[d]_{\vec d} \ar@<0.5ex>[d]^{\vec c} & & 
B_1 \ar@<-0.5ex>[d]_{d} \ar@<0.5ex>[d]^{c} \\
B_1 \ar[rr]_-{\delta_0=d} & & B_0} \;\;\;\;\;
\xymatrix{\vec B_1 \ar[rr]^-{\gamma_1=m_1 \cdot \pi_2} \ar@<-0.5ex>[d]_{\vec d} \ar@<0.5ex>[d]^{\vec c} & & 
B_1 \ar@<-0.5ex>[d]_{d} \ar@<0.5ex>[d]^{c} \\
B_1 \ar[rr]_-{\gamma_0=c} & & B_0}$$
Finally, the natural transformation $\beta \colon \delta \Rightarrow \gamma$ is simply $\beta = \id_{B_1} \colon B_1 \to B_1$.\\
To help intuition, let us point out that when $\cA$ is the category of sets, an element $S$ of the object $\vec B_1$ 
involved in the above description of the strong h-pullback is a commutative square
$$\xymatrix{\ar[r]^{g_0} \ar[d]_{b_1} & \ar[d]^{b_2} \\
\ar[r]_{f_0} \ar@{}[ru]|-{S} & }$$
with $m_1(S)=\langle b_1, f_0 \rangle, m_2(S) = \langle g_0, b_2 \rangle,
\vec d(S)=b_1, \vec c(S)=b_2, \delta_1(S)=g_0, \gamma_1(S)=f_0$.\\
Let us denote by $\tau_1 \colon \vec B_1 \rightarrow \vec B_1$ the unique morphism such that $\tau_1 \cdot m_1 = m_2$ and $\tau_1 \cdot m_2 = m_1$. Since $\tau_1 \cdot \tau_1 = \id_{\vec B_1}$, $\tau_1$ is an isomorphism. Together with $\tau_0 = \id_{B_1}$, we can transpose the groupoid structure of $\vec{\bB}$ to create the groupoid
$$\vec{\bB}'=(\xymatrix{\vec B_1 \times_{\vec{c} \,',\vec{d} \,'}\vec B_1 \ar[r]^-{\vec{m}'} & 
\vec B_1 \ar@<0.9ex>[r]^{\vec{d} \,'} \ar@<-0.9ex>[r]_{\vec{c} \,'} & B_1 \ar[l]|{\vec{e} \,'}} \;,\;
\xymatrix{\vec B_1 \ar[r]^{\vec{i} \,'} & \vec B_1})$$
and the isomorphism $\tau \colon \vec{\bB}' \rightarrow \vec{\bB}$. With this notation, we have $\vec{d} \,' = m_2 \cdot \pi_1$ and $\vec{c} \,' = m_1 \cdot \pi_2$.
We will need this isomorphism later on.
}\end{Text}

\begin{Text}\label{TextStrongHpbGrpdConstr}{\rm
Putting together~\ref{TextPBStrongHpb} and~\ref{TextStrongHpbGrpd}, we can conclude that the 2-category $\Grpd(\cA)$ 
has strong h-pullbacks. Moreover, a strong h-pullback
$$\xymatrix{\bP \ar[r]^{G'} \ar[d]_{F'} & \bA \ar[d]^{F} \\
\bC \ar@{}[ru]^{\varphi}|{\Rightarrow} \ar[r]_{G} & \bB}$$
of $F \colon \bA \to \bB$ and $G \colon \bC \to \bB$ in $\Grpd(\cA)$ is described by the following diagram in $\cA$, where the top and bottom faces are limit diagrams:
$$\xymatrix{& & P_1 \ar[rd]|-{\varphi_1} \ar[rrrd]^-{G_1'} 
\ar@<-0.5ex>[ddd]_>>>>>>>>>>{\underline d}|(.52)\hole \ar@<0.5ex>[ddd]^>>>>>>>>>>{\underline c}|(.5)\hole \ar[lld]_-{F_1'} \\
C_1 \ar[rd]^{G_1} \ar@<-0.5ex>[ddd]_{d} \ar@<0.5ex>[ddd]^{c} & & & \vec B_1 \ar[lld]_>>>>>>>>>{m_2 \cdot \pi_1} 
\ar@<-0.5ex>[ddd]_{m_1 \cdot \pi_1} \ar@<0.5ex>[ddd]^{m_2 \cdot \pi_2} \ar[rd]^{m_1 \cdot \pi_2} 
& & A_1 \ar[ld]_<<<<<<{F_1} \ar@<-0.5ex>[ddd]_{d}  \ar@<0.5ex>[ddd]^{c} \\
& B_1 \ar@<-0.5ex>[ddd]_<<<<<<<<{d} \ar@<0.5ex>[ddd]^<<<<<<<<{c} 
& & & B_1 \ar@<-0.5ex>[ddd]_<<<<<<<<<<{d} \ar@<0.5ex>[ddd]^<<<<<<<<<<{c} \\
& & P_0 \ar[lld]_<<<<<<{F_0'}|(.47)\hole|(.52)\hole \ar[rd]_{\varphi_0} \ar[rrrd]^-{G_0'}|(.31)\hole|(.35)\hole|(.65)\hole|(.69)\hole \\
C_0 \ar[rd]_{G_0} & & & B_1 \ar[lld]^{d} \ar[rd]_{c} & & A_0 \ar[ld]^{F_0} \\
& B_0 & & & B_0}$$
}\end{Text}

\begin{Text}\label{TextFFGrpd}{\rm
In $\Grpd(\cA)$, as in any 2-category, the notion of equivalence makes sense. Moreover, in $\Grpd(\cA)$ we have also available the notion 
of weak equivalence. From~\cite{BP, FractEV}, recall that a functor $F \colon \bA \to \bB$ between groupoids in $\cA$ is:
\begin{enumerate}
\item fully faithful if and only if the following diagram is a limit diagram;
$$\xymatrix{& & A_1 \ar[lld]_{d} \ar[d]^{F_1} \ar[rrd]^{c} \\
A_0 \ar[rd]_{F_0} & & B_1 \ar[ld]^{d} \ar[rd]_{c} & & A_0 \ar[ld]^{F_0} \\
& B_0 & & B_0}$$
\item an equivalence if it is fully faithful and, moreover, if the first row in one (equivalently, in both) of the following diagrams is a split 
epimorphism (the squares are pullbacks).
$$\xymatrix{A_0 \times_{F_0,d}B_1 \ar[r]^-{\beta_d} \ar[d]_{\alpha_d} & B_1 \ar[r]^-{c} \ar[d]^{d} & B_0 \\
A_0 \ar[r]_-{F_0} & B_0}
\;\;\;\;\;
\xymatrix{A_0 \times_{F_0,c}B_1 \ar[r]^-{\beta_c} \ar[d]_{\alpha_c} & B_1 \ar[r]^-{d} \ar[d]^{c} & B_0 \\
A_0 \ar[r]_-{F_0} & B_0}$$
The functors $\delta$ and $\gamma$ defined in~\ref{TextStrongHpbGrpd} are examples of equivalences.
\setcounter{functor}{\value{enumi}}
\end{enumerate}
If $\cA$ is a regular category, one can say that a functor $F \colon \bA \to \bB$ between groupoids in $\cA$ is:
\begin{enumerate}
\setcounter{enumi}{\value{functor}}
\item {\em essentially surjective} if $\beta_d \cdot c$ (equivalently $\beta_c \cdot d$) of the above diagrams is a regular epimorphism;
\item a {\em weak equivalence} if it is fully faithful, and essentially surjective.
\end{enumerate}
}\end{Text}

\section{Fibrations and strong h-pullbacks}\label{section fibrations}

In this section, we assume that the base category $\cA$ is regular. 

\begin{Text}\label{TextFibration}{\rm
Let us recall the terminology for fibrations (= opfibrations) between groupoids (compare with \cite[Definition 5.1]{EKVdL} for the notion of $\mathcal E$-fibrations between internal categories, w.r.t.\ a class $\mathcal E$ of morphisms of $\cA$).
Consider a functor $F \colon \bA \to \bB$ between groupoids in $\cA$, and the induced factorizations through the 
pullbacks as in the following diagrams.
$$\xymatrix{A_1 \ar[rr]^{F_1} \ar[rd]^{\tau_d}  \ar[dd]_{d} & & B_1 \ar[dd]^{d} \\
& A_0 \times_{F_0,d}B_1 \ar[ru]^{\beta_d} \ar[ld]_{\alpha_d} \\
A_0 \ar[rr]_{F_0} & & B_0}
\;\;\;\;\;\;\;
\xymatrix{A_1 \ar[rr]^{F_1} \ar[rd]^{\tau_c}  \ar[dd]_{c} & & B_1 \ar[dd]^{c} \\
& A_0 \times_{F_0,c}B_1 \ar[ru]^{\beta_c} \ar[ld]_{\alpha_c} \\
A_0 \ar[rr]_{F_0} & & B_0}$$
\begin{enumerate}
\item $F$ is a {\em fibration} when $\tau_d$ (equivalently, $\tau_c$) is a regular epimorphism.
\item $F$ is a {\em split epi fibration} when $\tau_d$ (equivalently, $\tau_c$) is a split epimorphism.
\item $F$ is a {\em discrete fibration} when $\tau_d$ (equivalently, $\tau_c$) is an isomorphism.
\end{enumerate}
}\end{Text}

In~\cite{Street}, Street defined $0$-fibrations in the more general context of a representable 2-category. 
As Chevalley criterion~\cite{Gray,Kock,Street}, he characterized $0$-fibrations as those $F$ for which the canonical functor 
$S \colon \vec \bA \to (F \downarrow \bB)$ (where $(F \downarrow \bB)$ is the comma object of $F$ over $\bB$, which coincides 
with the strong h-pullback of $F$ and $\Id_{\bB}$ in $\Grpd(\cA)$) has a left adjoint weak right inverse, i.e., the unit of the adjuntion is an isomorphism. 
One can show that a functor $F \colon \bA \to \bB$ between groupoids in $\cA$ is a split epi fibration if and only if the comparison functor $S$ has a left adjoint right inverse, i.e., the unit of the adjunction is an identity. Therefore, our notion of split epi fibration is a bit stronger than Street's notion of 
$0$-fibrations.

In the next characterization of fibrations and split epi fibrations, we use the canonical embedding $N \colon [B_0] \to \bB$ of the discrete groupoid 
of objects of $\bB$ into $\bB$. Explicitly, 
$$\xymatrix{B_0 \ar@<-0,5ex>[d]_{\id} \ar@<0,5ex>[d]^{\id} \ar[rr]^-{N_1=e} & & B_1 \ar@<-0,5ex>[d]_{d} \ar@<0,5ex>[d]^{c} \\
B_0 \ar[rr]_-{N_0=\id} & & B_0}$$
The discrete groupoid 2-functor is denoted as
$$[\ \ ]\colon \cA\to \Grpd(\cA).$$

\begin{Proposition}\label{PropCarattFibr}
Consider a functor $F \colon \bA \to \bB$ between groupoids in $\cA$, and the comparison functor $T$ from the pullback 
to the strong h-pullback, as in the following diagram.
$$\xymatrix{[B_0] \times_{N,F} \bA \ar@/^1pc/[rrrd]^{\widehat N} \ar[rd]^{T} \ar@/_1pc/[rdd]_{\widehat F} \\
& \bV(F) \ar[rr]^-{N'} \ar[d]_{F'} & \ar@{}[d]_{v(F)}^-{\Rightarrow} & \bA \ar[d]^{F} \\
& [B_0] \ar[rr]_-{N} & & \bB}$$
\begin{enumerate}
\item $F$ is a fibration if and only if $T$ is a weak equivalence.
\item $F$ is a split epi fibration if and only if $T$ is an equivalence.
\end{enumerate}
\end{Proposition}

In fact, we are going to prove a more precise statement: the arrow attesting that $T$ is essentially surjective is the same 
arrow $\tau_c$ attesting that $F$ is a fibration.

(To help intuition, it is worth providing a description of the groupoid $\mathbb{V}(F)$ when the base category $\cA$ is the category of sets. 
In this case, one has 
$$\mathbb{V}(F)=\underset{b\in\bB}{\coprod }(b\downarrow F)$$
i.e., the disjoint union of the comma categories $(b\downarrow F)$.)

\proof
Since pullbacks in $\Grpd(\cA)$ are strong (\ref{TextStrongHpbGrpd}), we can apply point~1 of~\ref{TextStability} 
and we know that $T$ is fully faithful. Now we have to compare $T_0' \cdot \underline c$ with $\tau_c$.
$$\xymatrix{A_0 \times_{T_0,\underline d} \bV(F)_1 \ar[d]_{\underline d'} \ar[r]^-{T_0'} & 
\bV(F)_1 \ar[d]^{\underline d} \ar[r]^-{\underline c} & \bV(F)_0 \\
A_0 \ar[r]_-{T_0} & \bV(F)_0}
\;\;\;\;\;
\xymatrix{ & A_1 \ar[ld]_{c} \ar[d]^{\tau_c} \ar[rd]^{F_1} \\
A_0 & A_0 \times_{F_0,c} B_1 \ar[l]^-{\alpha_c} \ar[r]_-{\beta_c} & B_1}$$
The diagram giving the strong h-pullback $\bV(F)$ is
$$
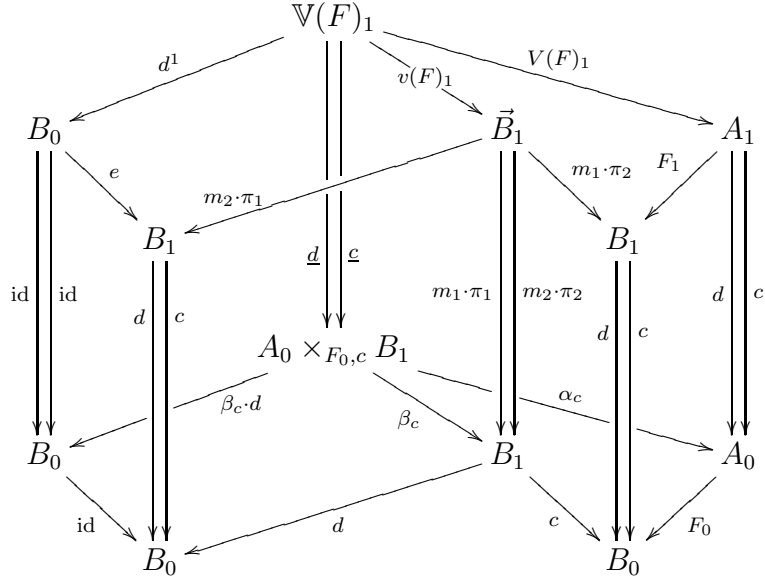
\xymatrix{& & \bV(F)_1 \ar[rd]|-{v(F)_1} \ar[rrrd]^-{V(F)_1} 
\ar@<-0.5ex>[ddd]_>>>>>>>>>>{\underline d}|(.51)\hole \ar@<0.5ex>[ddd]^>>>>>>>>>>{\underline c}|(.5)\hole \ar[lld]_-{d^1} \\
B_0 \ar[rd]^{e} \ar@<-0.5ex>[ddd]_{\id} \ar@<0.5ex>[ddd]^{\id} & & & \vec B_1 \ar[lld]_>>>>>>>>>{m_2 \cdot \pi_1} 
\ar@<-0.5ex>[ddd]_{m_1 \cdot \pi_1} \ar@<0.5ex>[ddd]^{m_2 \cdot \pi_2} \ar[rd]^{m_1 \cdot \pi_2} 
& & A_1 \ar[ld]_<<<<<<{F_1} \ar@<-0.5ex>[ddd]_{d}  \ar@<0.5ex>[ddd]^{c} \\
& B_1 \ar@<-0.5ex>[ddd]_<<<<<<<<{d} \ar@<0.5ex>[ddd]^<<<<<<<<{c} 
& & & B_1 \ar@<-0.5ex>[ddd]_<<<<<<<<<<{d} \ar@<0.5ex>[ddd]^<<<<<<<<<<{c} \\
& & A_0 \times_{F_0,c} B_1 \ar[lld]^<<<<<<{\beta_c \cdot d}|(.58)\hole|(.62)\hole \ar[rd]_{\beta_c} \ar[rrrd]^-{\alpha_c}|(.41)\hole|(.44)\hole|(.7)\hole|(.73)\hole \\
B_0 \ar[rd]_{\id} & & & B_1 \ar[lld]^{d} \ar[rd]_{c} & & A_0 \ar[ld]^{F_0} \\
& B_0 & & & B_0}$$
so that $T_0$ is the factorization through the pullback as in the following diagram.
$$
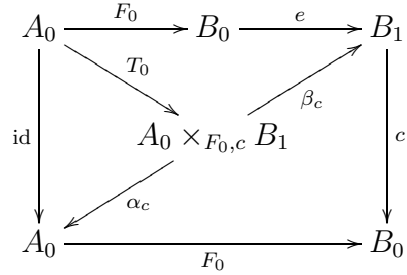
\xymatrix{A_0 \ar[dd]_{\id} \ar[rd]^{T_0} \ar[r]^-{F_0} & B_0 \ar[r]^-{e} & B_1 \ar[dd]^{c} \\
& A_0 \times_{F_0,c} B_1 \ar[ld]^{\alpha_c} \ar[ru]_{\beta_c}  \\
A_0 \ar[rr]_-{F_0} & & B_0}$$
Now we construct an arrow $\overline f \colon A_1 \to \bV(F)_1$ in three steps:
$$
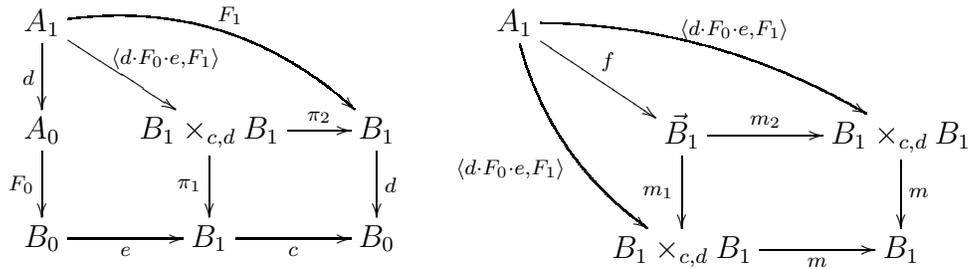
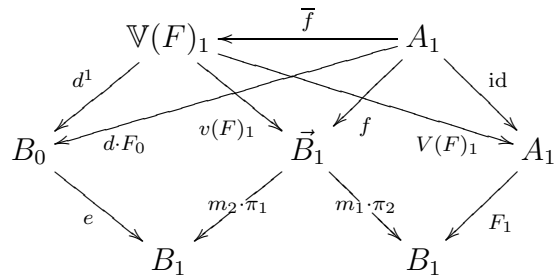
\xymatrix{A_1 \ar[d]_{d} \ar[rd]^{\langle d \cdot F_0 \cdot e,F_1 \rangle} \ar@/^1,5pc/[rrd]^{F_1} \\
A_0 \ar[d]_{F_0} & B_1 \times_{c,d} B_1 \ar[d]_{\pi_1} \ar[r]^-{\pi_2} & B_1 \ar[d]^{d} \\
B_0 \ar[r]_-{e} & B_1 \ar[r]_-{c} & B_0}
\;\;\;\;\;
\xymatrix{A_1 \ar@/_1pc/[rdd]_{\langle d \cdot F_0 \cdot e,F_1 \rangle} \ar[rd]^{f} \ar@/^1pc/[rrd]^{\langle d \cdot F_0 \cdot e,F_1 \rangle} \\
& \vec B_1 \ar[d]_{m_1} \ar[r]^-{m_2} & B_1 \times_{c,d} B_1 \ar[d]^{m} \\
& B_1 \times_{c,d} B_1 \ar[r]_-{m} & B_1}$$
$$\xymatrix{ & \bV(F)_1 \ar[ld]_{d^1} \ar[rd]_>>>>>{v(F)_1} \ar[rrrd]_>>>>>>>>{V(F)_1} & & 
A_1 \ar[ll]_-{\overline f} \ar[llld]^>>>>>>>>>{d \cdot F_0} \ar[ld]^>>>>{f} \ar[rd]^{\id} \\
B_0 \ar[rd]_{e} & & \vec B_1 \ar[ld]|{m_2 \cdot \pi_1} \ar[rd]|{m_1 \cdot \pi_2} & & A_1 \ar[ld]^{F_1} \\
& B_1 & & B_1}$$
Finally, we get the following diagram
$$\xymatrix{A_1 \ar[d]_{d} \ar[r]_-{\overline f} \ar@/^1,5pc/[rr]^{\tau_c} & \bV(F)_1 \ar[d]^{\underline d} \ar[r]_-{\underline c} & 
A_0 \times_{F_0,c} B_1 \\
A_0 \ar[r]_-{T_0} & A_0 \times_{F_0,c} B_1}$$
and we have to check that it is commutative and that the square is a pullback. Once this done, the commutativity of the upper region
immediately gives both statements of the proposition.\\
Commutativity of the upper region:
$$\overline f \cdot \underline c \cdot \alpha_c = \overline f \cdot V(F)_1 \cdot c = c = \tau_c \cdot \alpha_c$$
$$\overline f \cdot \underline c \cdot \beta_c = \overline f \cdot v(F)_1 \cdot m_2 \cdot \pi_2 = f \cdot m_2 \cdot \pi_2 = F_1 = \tau_c \cdot \beta_c$$
Commutativity of the square:
$$\overline f \cdot \underline d \cdot \alpha_c = \overline f \cdot V(F)_1 \cdot d = d = 
d \cdot T_0 \cdot \alpha_c$$
$$\overline f \cdot \underline d \cdot \beta_c = \overline f \cdot v(F)_1 \cdot m_1 \cdot \pi_1 = 
f \cdot m_1 \cdot \pi_1 = d \cdot F_0 \cdot e = d \cdot T_0 \cdot \beta_c$$
Universality of the square: consider the factorization of the square through the pullback
$$\xymatrix{A_1 \ar[dd]_{d} \ar[rd]^{\langle d, \overline f \rangle} \ar[rr]^-{\overline f} & & \bV(F)_1 \ar[dd]^{\underline d} \\
& P \ar[ld]^{p_1} \ar[ru]_{p_2} \\
A_0 \ar[rr]_-{T_0} & & A_0 \times_{F_0,c} B_1}$$
and the arrow
$$\xymatrix{P \ar[r]^-{p_2} & \bV(F)_1 \ar[r]^-{V(F)_1} & A_1}$$
Since $\overline f$ is a (split) monomorphism, in order to prove that $\langle d, \overline f \rangle$ and $p_2 \cdot V(F)_1$ realize 
an isomorphism, it is enough to check the conditions $p_2 \cdot V(F)_1 \cdot d = p_1$ and $p_2 \cdot V(F)_1 \cdot \overline f = p_2$. 
The first one is easy:
$$p_2 \cdot V(F)_1 \cdot d = p_2 \cdot \underline d \cdot \alpha_c = p_1 \cdot T_0 \cdot \alpha_c = p_1$$
For the second one, we compose with the three limit projections:
$$p_2 \cdot V(F) \cdot \overline f \cdot V(F)_1 = p_2 \cdot V(F)_1 \cdot \id = p_2 \cdot V(F)_1$$
$$p_2 \cdot d^1 = p_2 \cdot \underline d \cdot \beta_c \cdot d = p_1 \cdot T_0 \cdot \beta_c \cdot d = p_1 \cdot F_0 \cdot e \cdot d =
p_1 \cdot F_0 = p_1 \cdot T_0 \cdot \alpha_c \cdot F_0 =$$
$$= p_2 \cdot \underline d \cdot \alpha_c \cdot F_0 = 
p_2 \cdot V(F)_1 \cdot d \cdot F_0 = p_2 \cdot V(F)_1 \cdot \overline f \cdot d^1$$
and, when composing with $v(F)_1 \colon \bV(F)_1 \to \vec B_1$, we still have to compose with the four pullback projections 
out from $\vec B_1 \colon$
$$p_2 \cdot V(F)_1 \cdot \overline f \cdot v(F)_1 \cdot m_1 \cdot \pi_1 = p_2 \cdot V(F)_1 \cdot f \cdot m_1 \cdot \pi_1 =
p_2 \cdot V(F)_1 \cdot d \cdot F_0 \cdot e = p_2 \cdot \underline d \cdot \alpha_c \cdot F_0 \cdot e =$$
$$=p_1 \cdot T_0 \cdot \alpha_c \cdot F_0 \cdot e = p_1 \cdot F_0 \cdot e = p_1 \cdot T_0 \cdot \beta_c = 
p_2 \cdot \underline d \cdot \beta_c = p_2 \cdot v(F)_1 \cdot m_1 \cdot \pi_1$$
$$p_2 \cdot V(F)_1 \cdot \overline f \cdot v(F)_1 \cdot m_1 \cdot \pi_2 = p_2 \cdot V(F)_1 \cdot f \cdot m_1 \cdot \pi_2 = 
p_2 \cdot V(F)_1 \cdot F_1 = p_2 \cdot v(F)_1 \cdot m_1 \cdot \pi_2$$
$$p_2 \cdot V(F)_1 \cdot \overline f \cdot v(F)_1 \cdot m_2 \cdot \pi_2 = p_2 \cdot V(F)_1 \cdot f \cdot m_2 \cdot \pi_2 =
p_2 \cdot V(F)_1 \cdot F_1 =$$
$$=  p_2 \cdot \langle \underline d \cdot \alpha_c \cdot F_0 \cdot e, V(F)_1 \cdot F_1 \rangle \cdot m = (\star) =
p_2 \cdot \langle v(F)_1 \cdot m_1 \cdot \pi_1, v(F)_1 \cdot m_1 \cdot \pi_2 \rangle \cdot m =$$
$$=  p_2 \cdot v(F)_1 \cdot m_1 \cdot m =
p_2 \cdot v(F)_1 \cdot m_2 \cdot m = p_2 \cdot \langle v(F)_1 \cdot m_2 \cdot \pi_1, v(F)_1 \cdot m_2 \cdot \pi_2 \rangle \cdot m =$$
$$= p_2 \cdot \langle d^1 \cdot e, v(F)_1 \cdot m_2 \cdot \pi_2 \rangle \cdot m = p_2 \cdot v(F)_1 \cdot m_2 \cdot \pi_2$$
where in the $(\star)$-labelled step we use $p_2 \cdot \underline d \cdot \alpha_c \cdot F_0 \cdot e = p_2 \cdot v(F)_1 \cdot m_1 \cdot \pi_1$
from the first equality. As far as the last equality is concerned, observe that 
$$v(F)_1 \cdot m_2 \cdot \pi_2 \cdot d \cdot e = v(F)_1 \cdot m_2 \cdot \pi_1 \cdot c \cdot e =
d^1 \cdot e \cdot c \cdot e = d^1 \cdot e = v(F)_1 \cdot m_2 \cdot \pi_1.$$
Therefore, using the third equality, we have
$$p_2 \cdot V(F)_1 \cdot \overline f \cdot v(F)_1 \cdot m_2 \cdot \pi_1 = 
p_2 \cdot V(F)_1 \cdot \overline f \cdot v(F)_1 \cdot m_2 \cdot \pi_2 \cdot d \cdot e =$$
$$= p_2 \cdot v(F)_1 \cdot m_2 \cdot \pi_2 \cdot d \cdot e = p_2 \cdot v(F)_1 \cdot m_2 \cdot \pi_1$$
and the proof is complete.
\endproof

\section{$\ast$-Fibrations and strong h-kernels}\label{section *-fibrations}

In this section, we assume that the base category $\cA$ is regular and pointed. 

\begin{Text}\label{TextStrongHKer}{\rm
Since $\cA$ is pointed, as a special case of~\ref{TextStrongHpbGrpdConstr} we get a 
description of the strong h-kernel of a functor $F \colon \bA \to \bB$ between groupoids in $\cA$.
$$\xymatrix{\bK(F) \ar[rr]^{K(F)} \ar[d]_{0} & & \bA \ar[d]^{F} \\
[0] \ar@{}[rru]^{k(F)}|{\Rightarrow} \ar[rr]_{0} & & \bB}$$
Once again, let us make explicit the diagram in $\cA$ giving the strong h-kernel $\bK(F) \colon$
$$\xymatrix{& & \bK(F)_1 \ar[rd]|-{k(F)_1} \ar[rrrd]^-{K(F)_1} 
\ar@<-0.5ex>[ddd]_>>>>>>>>>>{\underline d}|(.51)\hole \ar@<0.5ex>[ddd]^>>>>>>>>>>{\underline c}|(.5)\hole \ar[lld]_-{0} \\
0 \ar[rd]^{0} \ar@<-0.5ex>[ddd]_{0} \ar@<0.5ex>[ddd]^{0} & & & \vec B_1 \ar[lld]_>>>>>>>>>{m_2 \cdot \pi_1} 
\ar@<-0.5ex>[ddd]_{m_1 \cdot \pi_1} \ar@<0.5ex>[ddd]^{m_2 \cdot \pi_2} \ar[rd]^{m_1 \cdot \pi_2} 
& & A_1 \ar[ld]_<<<<<<{F_1} \ar@<-0.5ex>[ddd]_{d}  \ar@<0.5ex>[ddd]^{c} \\
& B_1 \ar@<-0.5ex>[ddd]_<<<<<<<<{d} \ar@<0.5ex>[ddd]^<<<<<<<<{c} 
& & & B_1 \ar@<-0.5ex>[ddd]_<<<<<<<<<<{d} \ar@<0.5ex>[ddd]^<<<<<<<<<<{c} \\
& & \bK(F)_0 \ar[lld]_<<<<<<{0}|(.55)\hole|(.59)\hole \ar[rd]_{k(F)_0} \ar[rrrd]^(.5){K(F)_0}|(.36)\hole|(.39)\hole|(.67)\hole|(.7)\hole \\
0 \ar[rd]_{0} & & & B_1 \ar[lld]^{d} \ar[rd]_{c} & & A_0 \ar[ld]^{F_0} \\
& B_0 & & & B_0}$$
}\end{Text}

\begin{Text}\label{TextStarFibration}{\rm
Now we introduce $\ast$-fibrations and  split epi $\ast$-fibrations.
Consider a functor $F \colon \bA \to \bB$ between groupoids in $\cA$, and the induced factorizations 
$\widehat\tau_d$ and $\widehat\tau_c$ through the pullbacks
$$\xymatrix{\Ker(F_1 \cdot c) \ar[rr]^-{F_1^c} \ar[rd]^{\widehat \tau_d} \ar[d]_{k_{F_1 \cdot c}} 
& & \Ker(c) \ar[d]^{k_c} \\
A_1 \ar[d]_{d} & A_0 \times_{F_0,k_c \cdot d} \Ker(c) \ar[ru]^{\widehat \beta_d} \ar[ld]_{\widehat \alpha_d} 
& B_1 \ar[d]^{d} \\
A_0 \ar[rr]_-{F_0} & & B_0}$$
$$\xymatrix{\Ker(F_1 \cdot d) \ar[rr]^-{F_1^d} \ar[rd]^{\widehat \tau_c} \ar[d]_{k_{F_1 \cdot d}} 
& & \Ker(d) \ar[d]^{k_d} \\
A_1 \ar[d]_{c} & A_0 \times_{F_0,k_d \cdot c} \Ker(d) \ar[ru]^{\widehat \beta_c} \ar[ld]_{\widehat \alpha_c} 
& B_1 \ar[d]^{c} \\
A_0 \ar[rr]_-{F_0} & & B_0}$$
where $F_1^c$ and $F_1^d$ are determined by the conditions $F_1^c \cdot k_c = k_{F_1 \cdot c} \cdot F_1$
and $F_1^d \cdot k_d = k_{F_1 \cdot d} \cdot F_1$.
\begin{enumerate}
\item $F$ is a {\em $\ast$-fibration} when $\widehat\tau_d$ (equivalently, $\widehat\tau_c$) is a regular epimorphism.
\item $F$ is a {\em split epi $\ast$-fibration} when $\widehat\tau_d$ (equivalently, $\widehat\tau_c$) is a split epimorphism.
\end{enumerate}
}\end{Text}

\begin{Text}\label{TextDefStarFibr}{\rm
In order to justify the fact that in~\ref{TextStarFibration} one can equivalently use $\widehat\tau_d$ or $\widehat\tau_c$,
consider the isomorphisms $K(i)$ and $\overline{K(i)}$ obtained by
$$\xymatrix{\Ker(c) \ar[r]^-{k_c} \ar[d]_{K(i)} & B_1 \ar[r]^-{c} \ar[d]_{i} & B_0 \ar[d]^{\id} \\
\Ker(d) \ar[r]_-{k_d} & B_1 \ar[r]_-{d} & B_0}
\;\;\;\;\;\;
\xymatrix{\Ker(F_1 \cdot c) \ar[r]^-{k_{F_1 \cdot c}} \ar[d]_{\overline{K(i)}} & A_1 \ar[r]^-{F_1} \ar[d]_{i} & B_1 \ar[r]^-{c} \ar[d]^{i} & B_0 \ar[d]^{\id} \\
\Ker(F_1 \cdot d) \ar[r]_-{k_{F_1 \cdot d}} & A_1 \ar[r]_-{F_1} & B_1 \ar[r]_-{d} & B_0}$$
and then use them to build up the  commutative diagram
$$\xymatrix{\Ker(F_1 \cdot c) \ar[rr]^-{\widehat\tau_d} \ar[d]_{\overline{K(i)}} & & A_0 \times_{F_0,k_c \cdot d} \Ker(c) \ar[d]^{\id \times K(i)} \\
\Ker(F_1 \cdot d) \ar[rr]_-{\widehat\tau_c} & & A_0 \times_{F_0,k_d \cdot c} \Ker(d)}$$
}\end{Text}

\begin{Text}\label{TextFibrStarFibr}{\rm
Since in the diagram 
$$\xymatrix{\Ker(F_1 \cdot c) \ar[r]^-{\widehat\tau_d} \ar[d]_{k_{F_1 \cdot c}} & 
A_0 \times_{F_0,k_c \cdot d} \Ker(c) \ar[r]^-{\widehat\beta_d} \ar[d]_-{\id \times k_c} & \Ker(c) \ar[r]^-{0} \ar[d]_{k_c} & 0 \ar[d]^{0} \\
A_1 \ar@{}[ru]|{(1)} \ar[r]_-{\tau_d} & A_0 \times_{F_0,d} B_1 \ar@{}[ru]|{(2)} \ar[r]_-{\beta_d} & B_1 \ar@{}[ru]|{(3)} \ar[r]_-{c} & B_0}$$
part (2) and part (3) are pullbacks and the whole is a pullback (because $\tau_d \cdot \beta_d = F_1$), it follows
that part (1) also is a pullback. This proves that any fibration is a $\ast$-fibration and any split epi fibration is a split epi $\ast$-fibration.
}\end{Text}

\begin{Proposition}\label{PropCarattStarFibr}
Consider a functor $F \colon \bA \to \bB$ between groupoids in $\cA$, and the comparison $J$ from its kernel 
to its strong h-kernel as in the following diagram.
$$\xymatrix{\bKer(F) \ar@/^1pc/[rrrd]^{K_F} \ar[rd]^{J} \ar@/_1pc/[rdd]_{0} \\
& \bK(F) \ar[rr]^-{K(F)} \ar[d]_{0} & \ar@{}[d]_{k(F)}^{\Rightarrow} & \bA \ar[d]^{F} \\
& [0] \ar[rr]_-{0} & & \bB}$$
\begin{enumerate}
\item $F$ is a $\ast$-fibration if and only if $J$ is a weak equivalence.
\item $F$ is a split epi $\ast$-fibration if and only if $J$ is an equivalence.
\end{enumerate}
\end{Proposition}

Similarly to what we did in Proposition~\ref{PropCarattFibr}, we are going to prove a more precise statement: 
the arrow attesting that $J$ is essentially surjective is the same arrow $\widehat\tau_c$ attesting that $F$ is a $\ast$-fibration.

\proof
Since pullbacks in $\Grpd(\cA)$ are strong (\ref{TextStrongHpbGrpd}), we can apply point~1 of~\ref{TextStability} 
and we know that $J$ is fully faithful. Now we have to compare $J_0' \cdot \underline c$ with $\widehat\tau_c$.
$$\xymatrix{\Ker(F_0) \times_{J_0,\underline d} \bK(F)_1 \ar[d]_{\underline d'} \ar[r]^-{J_0'} & 
\bK(F)_1 \ar[d]^{\underline d} \ar[r]^-{\underline c} & \bK(F)_0 \\
\Ker(F_0) \ar[r]_-{J_0} & \bK(F)_0}
\;\;\;
\xymatrix{\Ker(F_1 \cdot d) \ar[rr]^-{\widehat\tau_c} & & A_0 \times_{F_0,k_d \cdot c} \Ker(d)}$$
From the explicit description of the strong h-kernel $\bK(F)$, we see that
$$\bK(F)_0 = A_0 \times_{F_0,k_d \cdot c} \Ker(d) \;,\;\; K(F)_0 = \widehat\alpha_c \;,\;\; k(F)_0 = \widehat\beta_c \cdot k_d$$
Moreover, the universal property of $\Ker(F_0)$ gives the following factorization $\overline d \colon$
$$\xymatrix{\Ker(F_1 \cdot d) \ar[r]^-{k_{F_1 \cdot d}} \ar[d]_{\overline d} & A_1 \ar[r]^-{F_1} \ar[d]^{d} & B_1 \ar[d]^{d} \\
\Ker(F_0) \ar[r]_-{k_{F_0}} & A_0 \ar[r]_-{F_0} & B_0}$$
Now we construct an arrow $\overline f \colon \Ker(F_1 \cdot d) \to \bK(F)_1$ in three steps
$$\xymatrix{\Ker(F_1 \cdot d) \ar[r]^-{k_{F_1 \cdot d}} \ar[rd]|{\langle 0,k_{F_1 \cdot d} \cdot F_1 \rangle} 
\ar@/_2pc/[rdd]_{0} & A_1 \ar@/^1pc/[rd]^{F_1} \\
& B_1 \times_{c,d} B_1 \ar[r]^-{\pi_2} \ar[d]_{\pi_1} & B_1 \ar[d]^{d} \\
& B_1 \ar[r]_-{c} & B_0}
\;\;\;\;\;
\xymatrix{\Ker(F_1 \cdot d) \ar@/_1pc/[rdd]_{\langle 0,k_{F_1 \cdot d} \cdot F_1 \rangle} \ar[rd]^{f} 
\ar@/^1pc/[rrd]^{\langle 0,k_{F_1 \cdot d} \cdot F_1 \rangle} \\
& \vec B_1 \ar[d]_{m_1} \ar[r]^-{m_2} & B_1 \times_{c,d} B_1 \ar[d]^{m} \\
& B_1 \times_{c,d} B_1 \ar[r]_-{m} & B_1}$$
$$\xymatrix{ & \bK(F)_1 \ar[ld]_{0} \ar[rd]_>>>>>{k(F)_1} \ar[rrrd]_>>>>>>>>>>>>{K(F)_1} & & 
\Ker(F_1 \cdot d) \ar[ll]_-{\overline f} \ar[llld]^>>>>>>>>>{0} \ar[ld]^>>>>>>{f} \ar[rd]^{k_{F_1 \cdot d}} \\
0 \ar[rd]_{0} & & \vec B_1 \ar[ld]|{m_2 \cdot \pi_1} \ar[rd]|{m_1 \cdot \pi_2} & & A_1 \ar[ld]^{F_1} \\
& B_1 & & B_1}$$
Finally, we get the following diagram
$$\xymatrix{\Ker(F_1 \cdot d) \ar[d]_{\overline d} \ar[r]_-{\overline f} \ar@/^1,5pc/[rr]^{\widehat\tau_c} 
& \bK(F)_1 \ar[d]^{\underline d} \ar[r]_-{\underline c} & A_0 \times_{F_0,k_d \cdot c} \Ker(d) \\
\Ker(F_0) \ar[r]_-{J_0} & A_0 \times_{F_0,k_d \cdot c} \Ker(d)}$$
and we have to check that it is commutative and that the square is a pullback. Once this done, the commutativity of the upper region
immediately gives both parts of the statement.\\
Commutativity of the upper region:
$$\overline f \cdot \underline c \cdot K(F)_0 = \overline f \cdot K(F)_1 \cdot c = k_{F_1 \cdot d} \cdot c = 
\widehat\tau_c \cdot \widehat\alpha_c = \widehat\tau_c \cdot K(F)_0$$
$$\overline f \cdot \underline c \cdot k(F)_0 = \overline f \cdot k(F)_1 \cdot m_2 \cdot \pi_2 = f \cdot m_2 \cdot \pi_2 =
k_{F_1 \cdot d} \cdot F_1 = F_1^d \cdot k_d = \widehat\tau_c \cdot \widehat\beta_c \cdot k_d = \widehat\tau_c \cdot k(F)_0$$
Commutativity of the square:
$$\overline f \cdot \underline d \cdot K(F)_0 = \overline f \cdot K(F)_1 \cdot d = k_{F_1 \cdot d} \cdot d = 
\overline d \cdot k_{F_0} = \overline d \cdot J_0 \cdot K(F)_0$$
$$\overline f \cdot \underline d \cdot k(F)_0 = \overline f \cdot k(F)_1 \cdot m_1 \cdot \pi_1 = 
f \cdot m_1 \cdot \pi_1 = 0 = \overline d \cdot 0 = \overline d \cdot J_0 \cdot k(F)_0$$
Universality of the square: consider the factorization of the square through the pullback
$$\xymatrix{\Ker(F_1 \cdot d) \ar[dd]_{\overline d} \ar[rd]^{\langle \overline d, \overline f \rangle} \ar[rr]^-{\overline f} 
& & \bK(F)_1 \ar[dd]^{\underline d} \\
& P \ar[ld]^{p_1} \ar[ru]_{p_2} \\
\Ker(F_0) \ar[rr]_-{J_0} & & \bK(F)_0}$$
Moreover, since
$$p_2 \cdot K(F)_1 \cdot F_1 \cdot d = p_2 \cdot K(F)_1 \cdot d \cdot F_0 = p_2 \cdot \underline d \cdot K(F)_0 \cdot F_0 =
p_1 \cdot J_0 \cdot K(F)_0 \cdot F_0 = p_1 \cdot k_{F_0} \cdot F_0 = 0$$
the universal property of $\Ker(F_1 \cdot d)$ gives the following factorization $\rho \colon$
$$\xymatrix{P \ar[r]^-{p_2} \ar[rrd]_{\rho} & \bK(F)_1 \ar[r]^-{K(F)_1} & A_1 \ar[r]^-{F_1} & B_1 \ar[r]^-{d} & B_0 \\
& & \Ker(F_1 \cdot d) \ar[u]_{k_{F_1 \cdot d}}}$$
Since $\overline f$ is a monomorphism (because $\overline f \cdot K(F)_1 = k_{F_1 \cdot d}$), in order to prove that 
$\langle \overline d, \overline f \rangle$ and $\rho$ realize an isomorphism, it is enough to check the conditions 
$\rho \cdot \overline d = p_1$ and $\rho \cdot \overline f = p_2$. 
The first one is easy, just compose with the monomorphism $k_{F_0} \colon$
$$\rho \cdot \overline d \cdot k_{F_0} = \rho \cdot k_{F_1 \cdot d} \cdot d = p_2 \cdot K(F)_1 \cdot d = 
p_2 \cdot \underline d \cdot K(F)_0 = p_1 \cdot J_0 \cdot K(F)_0 = p_1 \cdot k_{F_0}$$
For the second one, we compose with the limit projections $K(F)_1$ and $k(F)_1$ and, when composing with 
$k(F)_1 \colon \bK(F)_1 \to \vec B_1$, we still have to compose with the four pullback projections out from $\vec B_1 \colon$
$$\rho \cdot \overline f \cdot K(F)_1 = \rho \cdot k_{F_1 \cdot d} = p_2 \cdot K(F)_1$$
$$\rho \cdot \overline f \cdot k(F)_1 \cdot m_1 \cdot \pi_1 = \rho \cdot f \cdot m_1 \cdot \pi_1 = \rho \cdot 0 = 
p_1 \cdot 0 = p_1 \cdot J_0 \cdot k(F)_0 = p_2 \cdot \underline d \cdot k(F)_0 = p_2 \cdot k(F)_1 \cdot m_1 \cdot \pi_1$$
$$\rho \cdot \overline f \cdot k(F)_1 \cdot m_1 \cdot \pi_2 = \rho \cdot f \cdot m_1 \cdot \pi_2 = 
\rho \cdot k_{F_1 \cdot d} \cdot F_1 = p_2 \cdot K(F)_1 \cdot F_1 = p_2 \cdot k(F)_1 \cdot m_1 \cdot \pi_2$$
$$\rho \cdot \overline f \cdot k(F)_1 \cdot m_2 \cdot \pi_1 = \rho \cdot \overline f \cdot 0 = 
p_2 \cdot 0 = p_2 \cdot k(F)_1 \cdot m_2 \cdot \pi_1$$
$$\rho \cdot \overline f \cdot k(F)_1 \cdot m_2 \cdot \pi_2 = \rho \cdot f \cdot m_2 \cdot \pi_2 = 
\rho \cdot k_{F_1 \cdot d} \cdot F_1 = p_2 \cdot K(F)_1 \cdot F_1 =$$
$$= \langle 0, p_2 \cdot K(F)_1 \cdot F_1 \rangle \cdot m = p_2 \cdot \langle \underline d \cdot k(F)_0, K(F)_1 \cdot F_1 \rangle \cdot m =
p_2 \cdot \langle k(F)_1 \cdot m_1 \cdot \pi_1, k(F)_1 \cdot m_1 \cdot \pi_2 \rangle \cdot m =$$
$$= p_2 \cdot k(F)_1 \cdot m_1 \cdot m = p_2 \cdot k(F)_1 \cdot m_2 \cdot m = 
p_2 \cdot \langle k(F)_1 \cdot m_2 \cdot \pi_1, k(F)_1 \cdot m_2 \cdot \pi_2 \rangle \cdot m =$$
$$= p_2 \cdot \langle 0, k(F)_1 \cdot m_2 \cdot \pi_2 \rangle \cdot m = p_2 \cdot k(F)_1 \cdot m_2 \cdot \pi_2$$
\endproof

From~\ref{TextFibrStarFibr} and Proposition~\ref{PropCarattStarFibr}, we get the following result, 
announced in Proposition~4.4 of~\cite{SnailMMVShort}:

\begin{Corollary}\label{CorFibJWeakEq}
Let $F \colon \bA \to \bB$ be a fibration between internal groupoids. The canonical comparison 
$J \colon \bKer(F) \to \bK(F)$ from the kernel to the strong h-kernel is a weak equivalence.
If $F$ is a split epi fibration, then $J$ is an equivalence.
\end{Corollary}

\begin{Text}\label{TextBetterNonLinSnake}{\rm
Thanks to Proposition~\ref{PropCarattStarFibr}, we can slightly improve Proposition~4.6 in~\cite{SnailMMVShort}: 
assume that $\cA$ is a pointed regular category with reflexive coequalizers and consider a $\ast$-fibration 
$F \colon \bA \to \bB$ in $\Grpd(\cA)$, with $\bA, \bB$ and $\bK(F)$ proper (in Proposition~4.6 of~\cite{SnailMMVShort},
 $F$ is assumed to be a fibration and not just a $\ast$-fibration). There exists an exact sequence 
$$\pi_1(\bKer(F)) \to \pi_1(\bA) \to \pi_1(\bB) \to \pi_0(\bKer(F)) \to \pi_0(\bA) \to \pi_0(\bB)$$
where $\pi_1(\bA)$ is the internal group of automorphisms on the base point of $\bA$ and $\pi_0(\bA)$ is 
the object of connected components of $\bA$. (Here, the exactness at $B$ of 
$$\xymatrix{A \ar[r]^-{f} & B \ar[r]^-{g} & C}$$
means that $f$ factors as a regular epimorphism followed by the kernel of $g$).
Indeed, since $J \colon \bKer(F) \to \bK(F)$ is a weak equivalence, the arrows
$\pi_0(F) \colon \pi_0(\bKer(F)) \to \pi_0(\bK(F))$ and $\pi_1(F) \colon \pi_1(\bKer(F)) \to \pi_1(\bK(F))$ 
are isomorphisms (Lemma~4.5 in~\cite{SnailMMVShort}). Therefore, the above exact sequence immediately 
follows from the exact sequence 
$$\pi_1(\bK(F)) \to \pi_1(\bA) \to \pi_1(\bB) \to \pi_0(\bK(F)) \to \pi_0(\bA) \to \pi_0(\bB)$$
established in Section~3 of~\cite{SnailMMVShort}.
}\end{Text}

\begin{Text}\label{TextJEquiv}{\rm
Corollary~\ref{CorFibJWeakEq} can be obtained also from Proposition~\ref{PropCarattFibr}
without using the notion of $\ast$-fibration (\ref{TextFibrStarFibr} and \ref{PropCarattStarFibr}).
Indeed, consider the comparison functors $L$ and $I$ as in the following pullback and strong h-pullback diagrams.
$$\xymatrix{\bKer(F) \ar@/^1pc/[rrrd]^{K_F} \ar[rd]^{L} \ar@/_1pc/[rdd]_{0} \\
& [B_0] \times_{N,F} \bA \ar[rr]^-{\widehat N} \ar[d]_{\widehat F} & & \bA \ar[d]^{F} \\
& [B_0] \ar[rr]_-{N} & & \bB} \;\;\;
\xymatrix{\bK(F) \ar@/^1pc/[rrrd]^{K(F)} \ar[rd]^{I} \ar@/_1pc/[rdd]_{0} \\
& \bV(F) \ar[rr]^-{N'} \ar[d]_{F'} & \ar@{}[d]_{v(F)}^{\Rightarrow} & \bA \ar[d]^{F} \\
& [B_0] \ar[rr]_-{N} & & \bB}$$
We get a commutative diagram
$$\xymatrix{\bKer(F) \ar[r]^-{J} \ar[d]_{L} & \bK(F) \ar[d]^{I} \\
[B_0] \times_{N,F} \bA \ar[r]_-{T} & \bV(F)}$$
which in fact is a pullback and where $I \colon \bK(F) \to \bV(F)$ is a discrete fibration.
Now, the fact that $J$ is a weak equivalence if $F$ is a fibration (or an equivalence if $F$ is a split epi fibration)
follows from Proposition~\ref{PropCarattFibr} and the following general lemma on pullbacks in $\Grpd(\cA)$.
}\end{Text}

\begin{Lemma}\label{LemmaPBdiscrFibr}
Consider a pullback in $\Grpd(\cA)$
$$\xymatrix{\bX \ar[d]_{\widehat G} \ar[r]^-{\widehat F} & \bC \ar[d]^{G} \\
\bA \ar[r]_-{F} & \bB}$$
and assume that $G$ is a discrete fibration.
\begin{enumerate}
\item If $F$ is a weak equivalence, then $\widehat F$ is a weak equivalence.
\item If $F$ is an equivalence, then $\widehat F$ is an equivalence.
\end{enumerate}
\end{Lemma}

\proof
Since in $\Grpd(\cA)$ pullbacks are strong (\ref{TextStrongHpbGrpd}), we know from point~2 of~\ref{TextStability} that, if $F$ 
is fully faithful, so is $\widehat F$.\\
Now, the universal property of the pullback $A_0 \times_{F_0,d} B_1$ gives a unique arrow $\lambda$ making the following diagram 
commutative.
$$\xymatrix{X_0 \times_{\widehat{F}_0,d} C_1 \ar[d]_{\pi_d} \ar[rdd]^>>>>>>{\lambda} \ar[r]^-{\gamma_d} & C_1 \ar[d]_{d} \ar[rrdd]^{G_1} \\
X_0 \ar[rdd]_{\widehat{G}_0} \ar[r]^>>>>>>>>{\widehat{F}_0}|!{[u];[rd]}\hole & C_0 \ar[rrdd]^<<<<<<<<{G_0}|!{[d];[dr]}\hole \\
& A_0 \times_{F_0,d} B_1 \ar[d]^{\alpha_d} \ar[rr]^>>>>>>>>{\beta_d} & & B_1 \ar[d]^{d} \\
& A_0 \ar[rr]_-{F_0} & & B_0}$$
Consider the following commutative diagrams.
$$\xymatrix{X_0 \times_{\widehat{F}_0,d} C_1 \ar[d]_{\pi_d} \ar[r]^-{\gamma_d} & C_1 \ar[d]^{d} \\
X_0 \ar[d]_{\widehat{G}_0} \ar[r]^-{\widehat{F}_0} \ar@{}[ru]|{(2)} & C_0 \ar[d]^{G_0} \\
A_0 \ar[r]_-{F_0} \ar@{}[ru]|{(1)} & B_0}
\;\;\;\;\;\;
\xymatrix{X_0 \times_{\widehat{F}_0,d} C_1 \ar[d]_{\lambda} \ar[r]^-{\gamma_d} & C_1 \ar[d]^{G_1} \\
A_0 \times_{F_0,d} B_1 \ar[d]_{\alpha_d} \ar[r]^-{\beta_d} \ar@{}[ru]|{(4)} & B_1 \ar[d]^{d} \\
A_0 \ar[r]_-{F_0} \ar@{}[ru]|{(3)} & B_0}
\;\;\;\;\;\;
\xymatrix{X_0 \times_{\widehat{F}_0,d} C_1 \ar[d]_{\gamma_d} \ar[r]^-{\lambda} & A_0 \times_{F_0,d} B_1 \ar[d]^{\beta_d} \\
C_1  \ar[d]_{c} \ar[r]^-{G_1} \ar@{}[ru]|{(4)} & B_1 \ar[d]^{c} \\
C_0 \ar[r]_-{G_0} \ar@{}[ru]|{(5)} & B_0}$$
Since (1) and (2) are pullbacks, then (1)+(2) also is a pullback. But (1)+(2) = (3)+(4) and (3) is a pullback, so (4) also is a pullback.
Since (5) is a pullback (because $G$ is a discrete fibration), we conclude that (4)+(5) is a pullback. The proof is now obvious:\\
- If $\beta_d \cdot c$ is a regular epimorphism, so is $\gamma_d \cdot c$, and this proves part 1 of the statement.\\
- If $\beta_d \cdot c$ is a split epimorphism, so is $\gamma_d \cdot c$, and this proves part 2 of the statement.
\endproof

\section{Normalized fibrations and normalized $\ast$-fibrations}\label{section normalized fibrations}

\begin{Text}\label{TextCatNullHom}{\rm
From~\cite{GR}, recall that a {\em category with null-homotopies} \underline{$\cB$} is given by
\begin{itemize}
\item a category $\cB$,
\item for each morphism $f \colon A \rightarrow B$ in $\cB$, a set $\cH(f)$ (the set of null-homotopies on $f$),
\item for each triple of composable morphisms $f \colon A \rightarrow B$, $g \colon B \rightarrow C$, $h \colon C \rightarrow D$, a map
$$f \circ - \circ h \colon \cH(g) \rightarrow \cH(f \cdot g \cdot h), \; \mu \mapsto f \circ \mu \circ h.$$
\end{itemize}
(If $f = \id_B$ or $h= \id_C$, we write $\mu \circ h$ or $f \circ \mu$ instead of $f \circ \mu \circ h$.)\\
These data have to satisfy
\begin{enumerate}
\item the identity condition: given a morphism $f \colon A \rightarrow B$, one has $\id_A \circ \mu \circ \id_B=\mu$ for all $\mu \in \cH(f)$,
\item the associativity condition: given morphisms
$$\xymatrix{A' \ar[r]^-{f'} & A \ar[r]^-{f} & B \ar[r]^-{g} & C \ar[r]^-{h} & D \ar[r]^-{h'} & D'}$$
one has $(f' \cdot f) \circ \mu \circ (h \cdot h') = f' \circ (f \circ \mu \circ h) \circ h'$ for any $\mu \in \cH(g)$.
\end{enumerate}
}\end{Text}

\begin{Text}\label{TextStrongHKerArrExample}{\rm
For what concerns the present work, a relevant example of category with null-homotopies is the category ${\Grpd(\cA)}$ of internal groupoids 
in a pointed category $\cA$, with the natural transformations $0\Rightarrow F$ playing the role of null-homotopies.
}\end{Text}

\begin{Text}\label{TextStrongHKerArr}{\rm
The structure of category with null-homotopies is not rich enough to express the notion of strong h-pullback, but still, following
\cite{GR, SnailEV}, we can express the notion of strong homotopy kernel.
Let \underline{$\cB$} be a category with null-homotopies and let $f \colon A \rightarrow B$ be a morphism in $\cB$. A triple
$$(\Ker(f), \; \K(f) \colon \Ker(f) \rightarrow A, \; k(f) \in \cH(\K(f) \cdot f))$$
\begin{enumerate}
\item is a {\em homotopy kernel} (h-kernel, for short) of $f$ if for any triple $$(D, g \colon D \rightarrow A, \mu \in \cH(g \cdot f)),$$
there exists a unique morphism $g' \colon D \rightarrow \Ker(f)$ in \underline{$\cB$} such that $g' \cdot \K(f) = g$ and $g' \circ k(f) = \mu$,
\item is a {\em strong homotopy kernel} (strong h-kernel, for short) of $f$ if it is a  h-kernel of~$f$ and, moreover, 
for any triple $(D, h \colon D \rightarrow \Ker(f), \mu \in \cH(h \cdot \K(f)))$ such that $\mu \circ f = h \circ k(f)$, there exists a unique 
$\lambda \in \cH(h)$ such that $\lambda \circ \K(f) = \mu$.
\end{enumerate}
Notice that in~\cite{GR}, the identity condition in the definition of a category with null-\-ho\-mo\-to\-pies has been omitted. 
We think it should not, since it allows to prove that h-kernels and strong h-kernels are determined up to isomorphism by their universal properties.

Finally, let us remark that the definition of strong h-kernels given in~\ref{TextStrongHKer} is consistent with the one given here, applied 
to category with null-homotopies ${\Grpd(\cA)}$ for a finitely pointed complete category $\cA$.
}\end{Text}

\begin{Text}\label{TextArrA}{\rm
For a category $\cA$, we consider the arrow category $\Arr(\cA) \colon$ the objects are the arrows 
$a \colon A \to A_0$ in $\cA$ and the morphisms $(f,f_0) \colon a \to b$ are commutative squares of the form
$$\xymatrix{A \ar[r]^-{f} \ar[d]_{a} & B \ar[d]^{b} \\
A_0 \ar[r]_-{f_0} & B_0}$$
From~\cite{SnailEV}, recall that $\Arr(\cA)$ is a category with null-homotopies: a null-homotopy for an 
arrow $(f,f_0)$ is a diagonal, that is an arrow $d \colon A_0 \to B$ such that $a \cdot d = f$ and $d \cdot b = f_0$. 
If $\cA$ has finite limits and a zero object, then $\Arr(\cA)$ has kernels and strong h-kernels. The kernel of $(f,f_0)$ is just
the level-wise kernel.
$$\xymatrix{\Ker(f) \ar[r]^-{k_f} \ar[d]_{\K(a)} & A \ar[r]^-{f} \ar[d]_{a} & B \ar[d]^{b} \\
\Ker(f_0) \ar[r]_-{k_{f_0}} & A_0 \ar[r]_-{f_0} & B_0}$$
To construct the strong h-kernel of $(f,f_0)$, consider the factorization through the pullback
$$\xymatrix{A \ar[rr]^-{f} \ar[dd]_{a} \ar[rd]^{\partial(f,f_0)_0} & & B \ar[dd]^{b} \\
& A_0 \times_{f_0,b} B \ar[ld]^{b'} \ar[ru]_{f_0'} \\
A_0 \ar[rr]_-{f_0} & & B_0}$$
The strong h-kernel is then given by the triple
$$(\partial(f,f_0)_0\,,\quad (id,b')\,,\quad f_0')$$
conveniently described by the following diagram.
$$\xymatrix{A \ar[r]^-{\id} \ar[dd]_{\partial(f,f_0)_0} & A \ar[r]^-{f} \ar[dd]^(.6){a} & B \ar[dd]^{b} \\ \\
A_0 \times_{f_0,b} B \ar[r]_-{b'} \ar[rruu]^(.4){f_0'} & A_0 \ar[r]_-{f_0} & B_0}$$
}\end{Text}

\begin{Text}\label{TextNormFunct}{\rm
For a finitely complete pointed category $\cA$, the examples presented in~\ref{TextStrongHKerArrExample} and in~\ref{TextArrA} are related by a functor, called the \emph{normalization functor}
$$\cN \colon \Grpd(\cA) \to \Arr(\cA)$$
which sends an (internal) functor $F \colon \bA \to \bB$ to the commutative diagram
$$\xymatrix{\Ker(d) \ar[r]^{\K_d(F)} \ar[d]_{k_d} & \Ker(d) \ar[d]^{k_d} \\
A_1 \ar[d]_{c} & B_1 \ar[d]^{c} \\
A_0 \ar[r]_{F_0} & B_0}$$
where $\K_d(F)$ is defined by $\K_d(F) \cdot k_d = k_d \cdot F_1$.
Moreover, if $\alpha \colon 0 \Rightarrow F$ is a null-homotopy in $\Grpd(\cA)$, it gives rise to a null-homotopy $\cN(\alpha)$ of $\cN(F)$ in the following way.
The natural transformation $\alpha$ is represented by an arrow $\alpha \colon A_0 \rightarrow B_1$. Since in particular $\alpha \cdot d = 0$, it factorises as $\alpha = \cN(\alpha) \cdot k_d$. We thus already have $\cN(\alpha) \cdot k_d \cdot c = \alpha \cdot c = F_0$.
$$\xymatrix{\Ker(d) \ar[d]_-{k_d} \ar[r]^-{\K_d(F)} & \Ker(d) \ar[d]^-{k_d} \\ A_1 \ar[d]_-{c} & B_1 \ar[d]^-{c} \\ A_0 \ar[r]_-{F_0} \ar@{.>}[ruu]^(.6){\cN(\alpha)} \ar[ru]^(.6){\alpha} & B_0}$$
Let us prove that $k_d \cdot c \cdot \cN(\alpha) = \K_d(F)$. 
Since $k_d$ is a monomorphism, it is enough to show that $k_d \cdot c \cdot \alpha = k_d \cdot F_1$.
The naturality of $\alpha$ means that the square 
$$\xymatrix{A_1 \ar[r]^-{\langle d \cdot \alpha, F_1 \rangle} \ar[d]_-{\langle 0,c \cdot \alpha \rangle} & B_1 \times_{c,d} B_1 \ar[d]^-{m} \\ B_1 \times_{c,d} B_1 \ar[r]_-{m} & B_1}$$
commutes. Precomposing with $k_d$, this gives
$$k_d \cdot c \cdot \alpha = k_d \cdot c \cdot \alpha \cdot \langle d \cdot e, \id \rangle \cdot m
= k_d \cdot \langle 0, c \cdot \alpha \rangle \cdot m
= k_d \cdot \langle d \cdot \alpha, F_1 \rangle \cdot m=$$
$$= \langle 0, k_d \cdot F_1 \rangle \cdot m
= k_d \cdot F_1 \cdot \langle d \cdot e, \id \rangle \cdot m
= k_d \cdot F_1$$
as required. This mapping $\alpha \mapsto \cN(\alpha)$ is compatible with the action of morphisms on null-homotopies.
}\end{Text}

The following two lemmas will be useful later on.

\begin{Lemma}\label{lemma1}
Let $\cA$ be a finitely complete pointed category and $F\colon \bA\to\bB$ a fully faithful functor between internal groupoids in $\cA$. Then its normalization $\cN(F)$ is a pullback in~$\cA$.
\end{Lemma}

\proof
Let us consider the diagram  below, where all the squares pullbacks.
$$
\xymatrix@!=4.5ex{
&&\Ker(d)\ar[dl]\ar[dr]\\
&\Ker(d)\ar[ld]\ar[dr]\ar@{}[rr]|{(1)}&&A_1\ar[ld]\ar[dr]\\
0\ar[rd]
&&A_0\times_{F_0,d}B_1\ar[dl]\ar[dr]\ar@{}[rr]|{(2)}
&&B_1\times_{c,F_0}A_0\ar[ld]\ar[rd]\\
&A_0\ar[rd]_{F_0}
&&B_1\ar[dl]_{d}\ar[rd]^{c}\ar@{}[rr]|{(3)}
&&A_0\ar[ld]^{F_0}\\
&&B_0&&B_0
}
$$
Then, $\cN(F)$ a pullback since it is precisely the region $(1)+(2)+(3)$.
\endproof

\begin{Lemma}\label{lemma2}
Let $\cA$ be a finitely complete pointed category and $F\colon \bA\to\bB$ a functor between internal groupoids in $\cA$. Then its strong h-kernel $\bK(F)$ is a discrete fibration.
\end{Lemma}

\proof
Let us expand the prisma on the right of the diagram of~\ref{TextStrongHKer}.
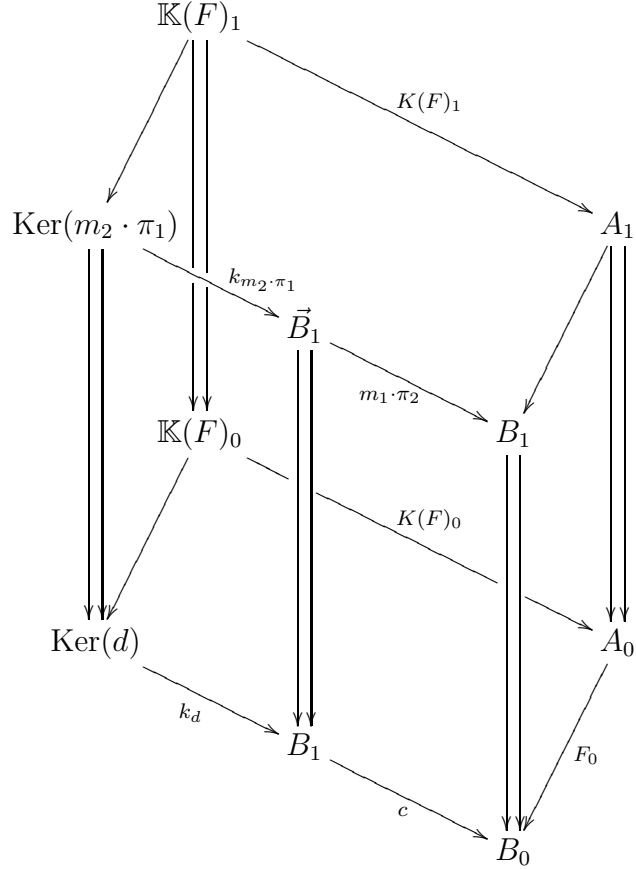
$$
\xymatrix@!=3ex{
&\bK(F)_1\ar@<-.5ex>[dddd]|(.62)\hole \ar@<+.5ex>[dddd]|(.63)\hole \ar[ddl]\ar[ddrrrr]^{K(F)_1}
\\
\\
\Ker(m_2\cdot\pi_1)\ar@<-.5ex>[dddd]\ar@<+.5ex>[dddd]\ar[drr]^(.7){k_{m_2\cdot\pi_1}}
&&&&&A_1\ar@<-.5ex>[dddd]\ar@<+.5ex>[dddd]\ar[ddl]
\\
&&\vec{B}_1\ar@<-.5ex>[dddd]\ar@<+.5ex>[dddd]\ar[drr]_{m_1\cdot\pi_2}
\\
&\bK(F)_0\ar[ddl]\ar[ddrrrr]^{K(F)_0}|(.23)\hole|(.26)\hole|(.73)\hole|(.76)\hole
&&&B_1\ar@<-.5ex>[dddd]\ar@<+.5ex>[dddd]
\\
\\
\Ker(d)\ar[drr]_{k_d}
&&&&&A_0\ar[ddl]^{F_0}
\\
&&B_1\ar[drr]_c
\\
&&&&B_0
}
$$
The upper and the lower faces are pullbacks by construction. 
The double front face over the arrow $k_d\cdot c$ can be interpreted as the normalization of the two functors
$$\xymatrix{\vec{\bB}' \ar@<-2pt>[d]_-{\tau \cdot \gamma}\ar@<2pt>[d]^-{\tau \cdot \delta} \\ \bB}$$
where the isomorphism $\tau \colon \vec{\bB}' \rightarrow \vec{\bB}$ have been described in~\ref{TextStrongHpbGrpd}.
These functors are equivalences, as pointed out in~\ref{TextFFGrpd}.
As a consequence, their normalizations are pullbacks in~$\cA$ by Lemma~\ref{lemma1}. Then, by elementary properties of pullbacks, also the double rear face over the arrow $K(F)_0$ must be made of two pullbacks, and this concludes our proof.
\endproof

\begin{Text}\label{TextPartialGrpd}{\rm
In order to define fully faithful morphisms in $\Arr(\cA)$, let us look more carefully to the situation in $\Grpd(\cA)$.
For a category $\cA$ with pullbacks and a functor $F \colon \bA \to \bB$ in $\Grpd(\cA)$, consider the following strong h-pullbacks 
$$\xymatrix{\vec{\bA} \ar[r]^{\gamma} \ar[d]_{\delta} & \bA \ar[d]^{\Id} \\
\bA \ar@{}[ru]^{\alpha}|{\Rightarrow} \ar[r]_{\Id} & \bA}
\;\;\;\;\;\;\;\;\;\;\;\;
\xymatrix{\bR(F) \ar[rr]^{\gamma(F)} \ar[d]_{\delta(F)} & & \bA \ar[d]^{F} \\
\bA \ar@{}[rru]^{\alpha(F)}|{\Rightarrow} \ar[rr]_{F} & & \bB}$$
and the unique functor $\partial(F) \colon \vec\bA \to \bR(F)$ such that
 $\partial(F) \cdot \delta(F) = \delta$, $\partial(F) \cdot \gamma(F) = \gamma$ and
$\partial(F) \cdot \alpha(F) = \alpha \cdot F$. The 0-level of the functor $\partial(F)$ is precisely the unique arrow
making commutative the following diagram.
$$\xymatrix{ & A_1 \ar[rr]^-{\partial(F)_0} \ar[ldd]_{d} \ar[rdd]_<<<<<<{F_1} \ar[rrrdd]_>>>>>>>>>>>>>>>>>>>>{c}
& & A_0 \times_{F_0,d} B_1 \times_{c,F_0} A_0 \ar[llldd]^>>>>>>>>>>>>>>{\delta(F)_0} \ar[ldd]^>>>>>>>>{\alpha(F)_0} \ar[rdd]^{\gamma(F)_0} \\ \\
A_0 \ar[rd]_{F_0} & & B_1 \ar[ld]^{d} \ar[rd]_{c} & & A_0 \ar[ld]^{F_0} \\
& B_0 & & B_0}$$
Therefore, we can say that the functor $F \colon \bA \to \bB$ is:
\begin{enumerate}
\item {\em faithful} if $\partial(F)_0$ is a monomorphism;
\item {\em full} if $\partial(F)_0$ is a regular epimorphism and in the context where $\cA$ is a regular category;
\item fully faithful if $\partial(F)_0$ is an isomorphism (accordingly to~\ref{TextFFGrpd}).
\end{enumerate}
}\end{Text}

\begin{Text}\label{TextPartialArr}{\rm
Now we imitate the previous argument using strong h-kernels in $\Arr(\cA)$. For a morphism $(f,f_0) \colon a \to b$ 
in $\Arr(\cA)$, consider the following strong h-kernels, together with the induced comparison arrow $\partial(f,f_0) \colon$
$$\xymatrix{\bK(\Id) \ar[r] \ar[d]_{\partial(f,f_0)} & a \ar[r]^-{\Id} \ar[d]_{\Id} & a \ar[d]^{(f,f_0)} \\
\bK(f,f_0) \ar[r] & a \ar[r]_-{(f,f_0)} & b}$$
The 0-level of $\partial(f,f_0)$ is precisely the factorization $\partial(f,f_0)_0 \colon A \to A_0 \times_{f_0,b}B$ 
through the pullback, as in~\ref{TextArrA}
(while the `domain level' of $\partial(f,f_0)$ is just the identity arrow on $A$). 
This suggests part of the following terminology.
}\end{Text}

\begin{Text} \label{TextClassesArr}{\rm
For a finitely complete pointed category $\cA$, consider an arrow $(f,f_0) \colon a \to b$ in $\Arr(\cA)$ together with the induced factorization 
$\partial(f,f_0)_0 \colon A \to A_0 \times_{f_0,b} B$ through the pullback, as in the description of the strong h-kernel in~\ref{TextArrA}. 
The arrow $(f,f_0)$ is:
\begin{enumerate}
\item {\em faithful} if $\partial(f,f_0)_0$ is a monomorphism;
\item {\em fully faithful} if $\partial(f,f_0)_0$ is an isomorphism;
\item {\em full} if $\partial(f,f_0)_0$ is a regular epimorphism and in the context where $\cA$ is a regular category;
\item {\em essentially surjective} if $f_0$ and $b$ are jointly strongly epimorphic;
\item a {\em weak equivalence} if it is fully faithful and essentially surjective;
\item a {\em fibration} if $f$ is a regular epimorphism and in the context where $\cA$ is a regular category.
\end{enumerate}
}\end{Text}

In order to compare the above terminology with the terminology for internal functors (\ref{TextFFGrpd}, \ref{TextFibration} 
and~\ref{TextPartialGrpd}), we need some intermediate steps. The first one is the version for strong h-pullbacks of the 
elementary fact that two parallel arrows in a pullback diagram have isomorphic kernels.

\begin{Lemma}\label{LemmaKerPB}
For a finitely complete pointed category $\cA$, consider the following diagram in $\Grpd(\cA)$, with the bottom square being a strong 
h-pullback, the right region being a strong h-kernel and the functor $\langle 0,K(F),k(F) \rangle$ determined by the universal property of the strong h-pullback, w.r.t.\ the triple $(0,K(F), k(F))$.
$$\xymatrix{\bK(F) \ar[r]^-{\Id} \ar[d]_{\langle 0,K(F),k(F) \rangle} & \bK(F) \ar[d]_{K(F)} \ar@/^4pc/[dd]^{0} \\
\bP \ar[d]_{F'} \ar[r]^-{G'}  & \bA \ar[d]_{F} \ar@{}[r]^{k(F)}|{\Leftarrow} &  \\
\bC \ar[r]_-{G} \ar@{}[ru]^{\varphi}|{\Rightarrow} & \bB}$$
Then the left column is a kernel.
\end{Lemma}

\proof
By the universal property of the strong h-pullback, we know that the functor $\langle 0,K(F),k(F) \rangle$ satisfies the conditions
$$
\langle 0,K(F),k(F) \rangle\cdot F'=0\,,\qquad
\langle 0,K(F),k(F) \rangle\cdot G' = K(F)\,,\qquad
\langle 0,K(F),k(F) \rangle\cdot \varphi=k(F).
$$
Let us re-write the diagram above in the following form.
$$
\xymatrix@C=12ex{
\bK(F)\ar[r]^{\langle 0,K(F),k(F) \rangle}\ar[d]&\bP \ar[d]_{F'} \ar[r]^-{G'}  & \bA \ar[d]_{F}   \\
0\ar[r]&\bC \ar[r]_-{G} \ar@{}[ru]^{\varphi}|{\Rightarrow} & \bB}
$$
By~\ref{TextPBStrongHpb} then, since the outer rectangle filled with the 2-cell $\langle 0,K(F),k(F) \rangle\cdot \varphi=k(F)$ is a strong h-pullback, then the left square is a strong pullback, i.e., $\langle 0,K(F),k(F) \rangle$ is the kernel of $F'$.
\endproof

\begin{Lemma}\label{LemmaPresStrongHKer}
If $\cA$ is a finitely complete pointed category, the normalization functor $\cN \colon \Grpd(\cA) \to \Arr(\cA)$ preserves kernels and strong h-kernels. 
\end{Lemma}

\proof
Preservation of kernels is an obvious argument of exchange of limits. Consider now a strong h-kernel in 
$\Grpd(\cA)$ and the canonical comparison $T=(t,t_0)$ with the strong h-kernel in $\Arr(\cA)$.
$$\xymatrix{\cN(\bK(F)) \ar[r]^-{\cN(K(F))} \ar[rd]_{T} & \cN(\bA) \ar[r]^-{\cN(F)} & \cN(\bB) \\
& \bK(\cN(F)) \ar[u] }$$
More explicitly, we get the following diagram in $\cA$
$$\xymatrix{\Ker(\underline d) \ar[d]_{k_{\underline d}} \ar[rd]_{t} \ar[rr]^-{t} & & 
\Ker(d) \ar[d]_{k_d} \ar[r]^-{\K_d(F)} & \Ker(d) \ar[d]^{k_d} \\
\bK(F)_1 \ar[d]_{\underline c} & \Ker(d) \ar[dd]^<<<<<<{\partial(\cN(F))_0} \ar[ru]_{\id} & A_1 \ar[d]_{c} & B_1 \ar[d]^{c} \\
\bK(F)_0 \ar[rd]_{t_0} \ar[rr]^(.3){K(F)_0}|!{[ur];[dr]}\hole & & A_0 \ar[r]_-{F_0} & B_0 \\
& A_0 \times_{F_0,k_d \cdot c} \Ker(d) \ar[ru]_{\beta'} }$$
where
$$\xymatrix{A_0 \times_{F_0,k_d \cdot c} \Ker(d) \ar[d]_{\beta'} \ar[r]^-{f_0'} & \Ker(d) \ar[d]^{k_d \cdot c} \\
A_0 \ar[r]_-{F_0} & B_0}$$
is a pullback, and $\partial(\cN(F))_0$ and $ t_0$ are determined by the conditions
$$\partial(\cN(F))_0 \cdot \beta' = k_d \cdot c \;,\;\; \partial(\cN(F))_0 \cdot f_0' = \K_d(F)$$
$$t_0 \cdot \beta' = K(F)_0 \;,\;\; t_0 \cdot f_0' \cdot k_d = k(F)_0 $$ 
and $t=\K_d(K(F))$.

Clearly, $t_0$ is an isomorphism (just look at the construction of $\bK(F)_0$ in~\ref{TextStrongHKer}). On the other hand, $t$ is the restriction to kernels of a pullback, since $K(F)$ is a discrete fibration by Lemma~\ref{lemma2}. Hence also $t$ is an isomorphism, and this concludes the proof.
\endproof

\begin{Lemma}\label{LemmaKerPBcat}
Let $\cA$ be a finitely complete pointed category and $F \colon \bA \to \bB$ be a functor in~$\Grpd(\cA)$. Consider the following diagram (notation as in~\ref{TextPartialGrpd}).
$$\xymatrix{\bK(\Id) \ar[rr]^-{\langle 0,K(\Id),k(\Id) \rangle} \ar[d]_{\langle 0,\K(\Id),k(\Id) \cdot F \rangle} 
& & \vec\bA \ar[r]^{\delta} \ar[d]_{\partial(F)} & \bA \ar[d]^{\Id} \\
\bK(F) \ar[rr]_-{\langle 0,K(F),k(F) \rangle} & & \bR(F) \ar[r]_{\delta(F)} & \bA}$$
Then the rows are kernels and the left-hand square is a pullback.
\end{Lemma}

\proof
Thanks to Lemma~\ref{LemmaPresStrongHKer}, $\langle 0,K(\Id),k(\Id) \cdot F \rangle \colon \bK(\Id) \rightarrow \bK(F)$ is sent by $\cN$ to
$$\partial(\cN(F)) \colon \cN(\bK(\Id))=\bK(\Id_{\cN(\bA)}) \rightarrow \bK(\cN(F)) = \cN(\bK(F)).$$
The fact that the rows are kernels is a particular case of 
Lemma~\ref{LemmaKerPB}. The last part of the statement now follows easily since $\Id \colon \bA \rightarrow \bA$ is a monomorphism.
\endproof

We are ready to compare the terminology in $\Grpd(\cA)$ and in $\Arr(\cA)$.
(See~\cite{BB} for the notions of protomodular and homological category. Compare also with~\cite{EKVdL} for points~9 and~10 
of the following result.)

\begin{Proposition}\label{PropFibrNorm}
Let $\cA$ be a finitely complete pointed category and $F \colon \bA \to \bB$ be a functor between groupoids in $\cA$.
\begin{enumerate}
\item If $F$ is faithful, then $\cN(F)$ is faithful.
\item If $F$ is fully faithful, then $\cN(F)$ is fully faithful.
\item If $\cA$ is regular and $F$ is full, then $\cN(F)$ is full.
\item If $\cA$ is protomodular and $\cN(F)$ is faithful, then $F$ is faithful.
\item If $\cA$ is protomodular and $\cN(F)$ is fully faithful, then $F$ is fully faithful.
\item If $\cA$ is regular and protomodular and $\cN(F)$ is full, then $F$ is full.
\item If $\cA$ is regular and protomodular and $F$ is essentially surjective, then $\cN(F)$ is essentially surjective.
\item If $\cA$ is regular and $\cN(F)$ is essentially surjective, then $F$ is essentially surjective.
\item If $\cA$ is regular and $F$ is a fibration, then $\cN(F)$ is a fibration.
\item If $\cA$ is regular and protomodular and $\cN(F)$ is a fibration, then $F$ is a fibration.
\end{enumerate}
\end{Proposition}

\proof
From 1 to 6. The 0-level of the diagram in Lemma~\ref{LemmaKerPBcat} gives the following diagram in $\cA \colon$
$$\xymatrix{\Ker(d) \ar[rr]^-{k_d} \ar[d]_{\partial(\cN(F))_0} & & A_1 \ar[rr]^{d} \ar[d]_{\partial(F)_0} & & A_0 \ar[d]^{\id} \\
A_0 \times_{F_0,k_d \cdot c} \Ker(d) \ar[rr]_-{\langle 0,K(F),k(f) \rangle_0} & & A_0 \times_{F_0,d} B_1 \times_{c,F_0} A_0 \ar[rr]_-{\delta(F)_0} & & A_0}$$
Thanks to Lemma~\ref{LemmaKerPBcat}, the rows are kernels and the left-hand square is a pullback. Therefore:\\
1. If $\partial(F)_0$ is a monomorphism, then $\partial(\cN(F))_0$ is also a monomorphism.\\
2. If $\partial(F)_0$ is an isomorphism, then $\partial(\cN(F))_0$ is also an isomorphism.\\
3. If $\partial(F)_0$ is a regular epimorphism, then $\partial(\cN(F))_0$ is also a regular epimorphism because the category $\cA$ is regular.\\
4. If $\partial(\cN(F))_0$ is a monomorphism, then $\partial(F)_0$ is also a monomorphism because in a protomodular category, pullbacks 
reflect monomorphisms (see Lemma~3.13 in~\cite{BG}).\\
5. If $\partial(\cN(F))_0$ is an isomorphism, then $\partial(F)_0$ is also an isomorphism because in a protomodular category, the 
Split Short Five Lemma holds (see Lemma~3.10 in~\cite{BG}).\\
6. If $\partial(\cN(F))_0$ is a regular epimorphism, then $\partial(F)_0$ is also a regular epimorphism by Proposition~8 in~\cite{DB} 
(see also Proposition~2.4 in~\cite{SnailEV}), which can be applied here because $\cA$ is regular and protomodular and 
$d \colon A_1 \to A_0$ is a split and then regular epimorphism.\\
7 and 8. Consider the following commutative diagram.
$$\xymatrix{A_0 \ar[rr]^-{\langle \id, F_0 \cdot e \rangle} \ar[rrd]_-{F_0} && A_0 \times_{F_0,d} B_1 \ar[d]^-{\beta_d \cdot c} &&
\Ker(d) \ar[ll]_-{\langle 0, k_d \rangle} \ar[lld]^-{k_d \cdot c} \\ && B_0 &&}$$
7. If $\beta_d \cdot c$ is a regular epimorphism, it is a strong one. Therefore, it suffices to prove that $\langle \id, F_0 \cdot e \rangle$ and
$\langle 0, k_d \rangle$ are jointly strongly epimorphic. This is the case since $\cA$ is protomodular and $\langle 0, k_d \rangle$ and $\langle \id, F_0 \cdot e \rangle$
are respectively a kernel and a section of $\alpha_d$ (see~\cite{DB,BG}).
$$\xymatrix{\Ker(d) \ar[r]^-{\langle 0, k_d \rangle} \ar[d] & A_0 \times_{F_0,d} B_1 \ar[d]_-{\alpha_d} \\ 0 \ar[r] 
& A_0 \ar@<-2pt>@/_/[u]_-{\langle \id, F_0 \cdot e \rangle}}$$
8. If $F_0$ and $k_d \cdot c$ are jointly strongly epimorphic, $\beta_d \cdot c$ is a strong epimorphism and so a regular epimorphism since $\cA$ is regular.\\
9 and 10. Consider the commutative diagram
$$\xymatrix{\Ker(d) \ar[r]^-{k_d} \ar[d]_{\K_d(F)} & A_1 \ar[d]_{\tau_d} \ar[r]^-{d} & A_0 \ar[d]^{\id} \\
\Ker(d) \ar[r]_-{\langle 0, k_d \rangle} & A_0 \times_{F_0,d}B_1 \ar[r]_-{\alpha_d} & A_0}$$
Since $\id \colon A_0 \to A_0$ is a monomorphism, the left-hand square is a pullback. Therefore:\\
9. If $\tau_d$ is a regular epimorphism, then $\K_d(F)$ is also a regular epimorphism because the category $\cA$ is regular.\\
10. If $\K_d(F)$ is a regular epimorphism, then $\tau_d$ is also a regular epimorphism by Proposition~8 in~\cite{DB}.
\endproof

We have not yet discussed $\ast$-fibrations in $\Arr(\cA)$. For this, we need a last preparatory step.
Given a morphism $(f,f_0) \colon a \to b$ in $\Arr(\cA)$, the triple $$(\bKer(f,f_0),\quad k_{(f,f_0)},\quad 0 \colon \Ker(f_0) \rightarrow B)$$ induces a canonical comparison  $J \colon \bKer(f,f_0) \to \bK(f,f_0)$ through the strong h-kernel of~$(f,f_0)$.

\begin{Lemma}\label{LemmaStrongKer}
In the category with null-homotopies $\Arr(\cA)$ for a finitely complete poin\-ted category $\cA$, kernels are strong, i.e., for any morphism $(f,f_0) \colon a \to b$,
the canonical comparison $J \colon \bKer(f,f_0) \to \bK(f,f_0)$ is fully faithful (see point~1 of~\ref{TextStability}).
\end{Lemma}

\proof
Using the descriptions of the kernel and of the strong h-kernel given in~\ref{TextArrA}, the comparison $J$ turns out to be
the left-hand square in the following diagram.
$$\xymatrix{\Ker(f) \ar[r]^-{k_f} \ar[d]_{\K(a)} & A \ar[d]_{\partial(f,f_0)_0} \ar[r]^-{f} & B \ar[d]^{\id} \\
\Ker(f_0) \ar[r]_-{\langle k_{f_0},0 \rangle} & A_0 \times_{f_0,b}B \ar[r]_-{f_0'} & B}$$
Since both rows are kernels and $\id \colon B \to B$ is a monomorphism, the left-hand square is a pullback, which means 
that $J$ is fully faithful.
\endproof

\begin{Text}\label{TextStarFibrArr}{\rm
Having in mind Proposition~\ref{PropCarattStarFibr}, we could now define a {\em $\ast$-fibration} in $\Arr(\cA)$ as a morphism 
$(f,f_0) \colon a \to b$ such that the canonical comparison $J \colon \bKer(f,f_0) \to \bK(f,f_0)$ is essentially surjective 
(and then, by Lemma~\ref{LemmaStrongKer}, a weak equivalence). Thanks to Proposition~\ref{PropFibrNorm}, if $\cA$ is homological, with such a 
notion of $\ast$-fibration we have that a functor $F$ is a $\ast$-fibration in $\Grpd(\cA)$ if and only if the morphism $\cN(F)$ is a 
$\ast$-fibration in~$\Arr(\cA)$. \\
Now that we have of the notions of fibration and $\ast$-fibration available in $\Arr(\cA)$, we can ask if every 
fibration is a $\ast$-fibration (this is the case in $\Grpd(\cA)$, as observed in~\ref{TextFibrStarFibr}). 
The surprise is that not only the answer is negative, but the expected implication ``fibration $\Rightarrow$ $\ast$-fibration'' 
is in fact equivalent to the condition of protomodularity.
}\end{Text}

\begin{Proposition}\label{PropCharProtomod}
The following conditions on a pointed regular category $\cA$ are equivalent:
\begin{enumerate}
\item $\cA$ is protomodular (and then homological).
\item For every fibration $(f,f_0) \colon a \to b$ in $\Arr(\cA)$, the canonical comparison 
$$J \colon \bKer(f,f_0) \to \bK(f,f_0)$$ is a weak equivalence.
\end{enumerate}
\end{Proposition}

\proof
1 $\Rightarrow$ 2. Suppose $\cA$ is homological and $(f,f_0)$ is a fibration in $\Arr(\cA)$.
Thanks to Lemma~\ref{LemmaStrongKer}, we already know that $J$ is full and faithful.
Consider now the following diagram where $\id \times f$ is a regular epimorphism since so is $f$.
$$\xymatrix{\Ker(f_0) \ar[rr]^-{\langle k_{f_0}, 0 \rangle} \ar[rrd]_-{\langle k_{f_0}, 0 \rangle} && A_0 \times_{f_0,a \cdot f_0} A \ar@{->>}[d]^-{\id \times f} &&
A \ar[ll]_-{\langle a, \id \rangle} \ar[lld]^-{\partial(f,f_0)_0} \\ && A_0 \times_{f_0,b} B &&}$$
Thus, in order to prove that $J$ is essentially surjective, it suffices to notice that the protomodularity of $\cA$ implies that
$\langle k_{f_0}, 0 \rangle$ and $\langle a, \id \rangle$ are jointly strongly epimorphic since they are respectively the kernel and a section of
the second projection $A_0 \times_{f_0,a \cdot f_0} A \rightarrow A$. \\
2 $\Rightarrow$ 1. Firstly, let us prove that if the kernel of a morphism $f_0 \colon A_0 \rightarrow B_0$ is the zero object, then $f_0$ is a monomorphism.
In order to do so, consider the fibration $(\id,f_0) \colon \id \rightarrow f_0$ in~$\Arr(\cA)$.
$$\xymatrix{A_0 \ar[r]^-{\id} \ar[d]_-{\id} & A_0 \ar[d]^-{f_0} \\ A_0 \ar[r]_-{f_0} & B_0 } \qquad
\xymatrix{0 \ar[r] \ar[d] & A_0 \ar[d]^-{\langle \id, \id \rangle} \\ 0 \ar[r] & A_0 \times_{f_0,f_0} A_0}$$
The diagram on the right represents the canonical comparison $J\colon \bKer(\id,f_0) \to \bK(\id,f_0)$.
By the assumption, we know that $0 \rightarrow A_0 \times_{f_0,f_0} A_0$ and $\langle \id, \id \rangle$ are jointly strongly epimorphic.
This is equivalent to the fact that $\langle \id, \id \rangle$ is a regular epimorphism. Since it is also a split monomorphism,
it is an isomorphism, which means that $f_0$ is a monomorphism.
Let us now prove that the Short Five Lemma holds in $\cA$. Consider the following diagram where both rows are kernel of regular epimorphisms and $\K(a)$ and $b$ are isomorphisms.
$$\xymatrix{\Ker(f) \ar[r]^-{k_f} \ar[d]_-{\K(a)} & A \ar@{->>}[r]^-{f} \ar[d]^-{a} & B \ar[d]^-{b} \\ \Ker(f_0) \ar[r]_-{k_{f_0}} & A_0 \ar@{->>}[r]_-{f_0} & B_0}$$

Since $\K(a)$ is a monomorphism, its kernel is the zero object. Since $b$ is a monomorphism, the left-hand square is a pullback, hence also the kernel of $a$ is zero. But, by the first part of the proof, this means that $a$ is a monomorphism. So, it remains to prove that it is a regular epimorphism.
The morphism $(f,f_0) \colon a \rightarrow b$ is a fibration in $\Arr(\cA)$. Thus, the comparison morphism $J \colon \bKer(f,f_0) \to \bK(f,f_0)$ is a weak equivalence.
Since $b$ is an isomorphism, this implies that $k_{f_0}$ and $a$ are jointly strongly epimorphic.
But since $\K(a)$ is an isomorphism, $k_{f_0}$ factors through $a$, so that $a$ is a regular epimorphism.
\endproof

\end{document}